\documentclass[12pt]{amsart}

\usepackage{amssymb}
\usepackage[all]{xy}

\newtheorem{propo}{Proposition}[section]
\newtheorem{defi}[propo]{Definition}

\newtheorem{lemma}[propo]{Lemma}
\newtheorem{corol}[propo]{Corollary}
\newtheorem{theor}[propo]{Theorem}
\newtheorem{theo}[propo]{Theorem}

\newtheorem{remar}[propo]{Remark}

\newcommand{\bl}{\begin{lemma}}
\newcommand{\el}{\end{lemma}}

\newcommand{\ra}{ \rightarrow }

\newcommand{\ld}{, \ldots ,}

\newcommand{\ov}{\overline}

\newcommand{\al}{\alpha}

\newcommand{\ep}{\varepsilon}

\newcommand{\lam}{\lambda }

\newcommand{\om}{\omega }

\newcommand{\diag}{\mathop{\rm diag}\nolimits}
\newcommand{\Hom}{\mathop{\rm Hom}\nolimits}
\newcommand{\Aut}{\mathop{\rm Aut}\nolimits}
\newcommand{\Inn}{\mathop{\rm Inn}\nolimits}
\newcommand{\Ker}{\mathop{\rm Ker}\nolimits}

\newcommand{\Irr}{\mathop{\rm Irr}\nolimits}

\newcommand{\St}{\mathop{\rm St}\nolimits}
\newcommand{\GL}{\mathop{\rm GL}\nolimits}
\newcommand{\PGL}{\mathop{\rm PGL}\nolimits}
\newcommand{\SL}{\mathop{\rm SL}\nolimits}
\newcommand{\SO}{\mathop{\rm SO}\nolimits}
\newcommand{\Sp}{\mathop{\rm Sp}\nolimits}
\newcommand{\ESp}{\mathop{\rm ESp}\nolimits}

\newcommand{\infl}{\mathop{\rm Infl}\nolimits}

\newenvironment{prf}{{\bf Proof.}}{\hfill $\Box$ \medskip}

\begin{document}

\title{The Weil-Steinberg character of finite classical groups}

\author{G.~Hiss and A.~Zalesski}

\address{G.H.: Lehrstuhl D f{\"u}r Mathematik, RWTH Aachen University,
52056 Aachen, Germany}
\address{A.Z.: School of Mathematics, University of East Anglia, Norwich,
NR47TJ, UK}

\email{G.H.: gerhard.hiss@math.rwth-aachen.de}
\email{A.Z.: alexandre@azalesski.wanadoo.co.uk}


\begin{abstract}
We compute the irreducible constitutents of the product of the Weil character
and the Steinberg character
in those finite classical groups for which a Weil character is defined,
namely the symplectic, unitary and general linear groups.
It turns out that this product is multiplicity free for the symplectic
and general unitary groups, but not for the general linear groups.

As an application we show that the restriction of the Steinberg character
of such a group to the subgroup stabilizing a vector in the natural module
is multiplicity free. The proof of this result for the unitary groups
uses an observation of Brunat, published as an appendix to our paper.

As our ``Weil character'' for the symplectic groups in even characteristic
we use the $2$-modular Brauer character of the generalized spinor representation.
Its product with the Steinberg character is the Brauer character of a projective
module. We also determine its indecomposable direct summands.
\end{abstract}

\maketitle

\section{Introduction}

\markboth{WEIL TIMES STEINBERG}{GERHARD HISS AND ALEX ZALESSKI}


The Steinberg character of a finite group of Lie type plays a prominent
role in its representation theory. During the recent two decades numerous 
papers have proved the significance of the Weil characters, although these 
are defined only for classical groups.

In this paper we study the product of the Weil characters with the Steinberg
character. For brevity we refer to such a procuct as the Weil-Steinberg 
character. Our main result claims that the decomposition of the Weil-Steinberg
character as sum of ordinary irreducible characters is multiplicity free for 
the symplectic and the unitary groups. In fact we provide a lot of information
about these irreducible constituents. One of the striking consequence is that the
Weil-Steinberg character is very much similar to the Gelfand-Graev character,
in the sense that the majority of the irreducible constituents of the latter
occur in the former and conversely.

Thus the Weil-Steinberg character can be viewed as a kind of deformation of 
the Gelfand-Graev character. (However, we do not think that the method used
for proving that the Gelfand-Graev character is multiplicity free can be used
for proving our result for the Weil-Steinberg character.) As the
Gelfand-Graev character plays a fundamental role in the representation theory
of groups of Lie type, one could expect that the Weil-Steinberg character will 
also appear significant.

At the moment we have two applications of our results. The first one is on
the restriction of the Steinberg character to the stabilizer of a vector of 
the natural module. We deduce that this restriction is multiplicity free. In 
addition we provide significant 
information on its irreducible constituents. We hope that this will stimulate
progress in the long-standing open problem of computing the restriction of an
arbitrary representation to the parabolic subgroup that is the stabilizer of 
an isotropic line of the natural module. Note that our proof used substantially
the ideas of the work of Jianbei An and the first author \cite{AH}, who obtain
this result for small-dimensional symplectic groups.

The second line of application of our method could be to computing 
decomposition numbers. The Weil-Steinberg character is the character
of the lift of a projective module in the defining characterisitic,
which is the direct sum of some principal indecomposable modules (PIMs for 
brevity). A straightforward consequence of our results is that each of these 
PIMs decomposes multiplicity freely as sum of ordinary irreducible characters,
hence certain columns of the decomposition matrix consist of the numbers~$1$ 
and~$0$ only. We do not determine these PIMs here but there are hints
that the number of them is not too small.  

Formally the Weil character cannot be defined for symplectic groups in 
characteristic~$2$. However the Brauer character of a certain module (which we
call the generalized spinor module) is an analogue of the Weil character in odd 
characteristic. Using this analogy, we obtain a similar result for symplectic
groups in even characteristic, namely, we show that the product of the 
generalized spinor Brauer character with the Steinberg character is 
multiplicity free when decomposed as sum of ordinary irreducible 
characters. In contrast with the odd characterisitic case, we also decompose 
this product as a direct sum of PIMs.

Before we state our main result, we need to specify precisely what we mean
by the Weil character in each case.

\begin{defi}\label{DefinitionOmega}
{\rm
Let $n > 1$ be an integer,~$q$ a power of the prime~$p$ and let $G = G_n(q)$
denote one of the following groups: $\Sp(2n,q)$, $U(2n,q)$, $U(2n+1,q)$, or
$\GL(n,q)$. 

(1) If $G = \Sp(2n,q)$ with $q$ odd we let~$\omega$ denote the character of 
one of (the two) Weil representation of~$G$ as introduced by 
G{\'e}rardin~\cite{G}.

(2) If $G = Sp(2n,q)$ with $q$ even we let~$\omega$ denote the class function
obtained be extending the Brauer character of the generalized spinor 
representation $\sigma_n$ of~$G$ by zeros on all of~$G$. (For a precise 
definition see Subsection~\ref{ModularCase} below.)

(3) If $G$ is a unitary group we let~$\hat{\omega}$ denote the character of 
the (unique) Weil representation of~$G$ as introduced by G{\'e}rardin~\cite{G}, 
and define $\omega$ by $\omega := \hat{\omega}$ if~$q$ is even, and by 
$\omega(g) := \det(g)^{(q+1)/2}\hat{\omega}(g)$, $g \in G$ if~$q$ is odd.

(4) If $G = \GL(n,q)$ we let $\hat{\omega}$ denote the permutation character
of~$G$ on its natural module, and define $\omega$ by $\omega := \hat{\omega}$
if~$q$ is even, and by $\omega(g) := \det(g)^{(q-1)/2}\hat{\omega}(g)$, $g \in
G$ if~$q$ is odd.
}
\end{defi}

\noindent 
In each case, $\omega$ is a class function of~$G$ of degree $q^n$, in fact
$\omega$ is a character of~$G$ except in Case~(2). We are interested in the
product $\omega \cdot \St$, where~$\St$ denotes the Steinberg character of~$G$.
Since the Steinberg character vanishes on $p$-singular elements, only the 
values of $\omega$ on $p$-regular, i.e., semisimple elements of~$G$ are 
relevant. (The two Weil characters of a symplectic group in odd characteristic
have the same restriction to the set of semisimple elements, so our choice
made in Case~(1) of Definition~\ref{DefinitionOmega} is not effective.)
Let~$V$ be the natural module for~$G$, and let $g \in G$. Write 
$N(V;g) := \dim\,\Ker( g - 1 )$ for the dimension of the $1$-eigenspace
of~$g$ on~$V$. Then if $g \in G$ is semisimple, we have $\omega(g) = 
\pm q^{N(V;g)/2}$ if~$G$ is a symplectic group, and $\omega(g) = 
\pm q^{N(V;g)}$, otherwise.
(For the sign in the Cases~(1) and~(3) of Definition~\ref{DefinitionOmega} 
see \cite[Corollaries 4.8.1, 4.8.2]{G}.)

The product $\omega \cdot \St$ is an ordinary character of~$G$, even in Case~(2)
of Definition~\ref{DefinitionOmega}. Since $\St$ is of $p$-defect~$0$, its
product with any ordinary character or (extended) $p$-modular character as 
in Case~(2) is the character of the lift of a projective module of~$G$ in
characteristic~$p$.

We can now formulate the main result of our paper.

\begin{theor}\label{Thm1Intro}
Let~$q$ be a power of the prime~$p$.
For a non-negative integer~$m$ let $G_m(q)$ denote one of the following groups:
$\Sp(2m,q)$, $U(2m,q)$, $U(2m+1,q)$, or $\GL(m,q)$ (with the convention that
$G_0(q)$ is the trivial group).

Fix a positive integer~$n > 1$, put $G := G_n$, and denote by~$V$ the natural
module for~$G$. Let $P_m$ denote the stabilizer in~$G$ of a totally
isotropic subspace of $V$ of dimension~$m$, so that the Levi subgroup
of~$P_m$ equals $\GL(m,q) \times G_{n-m}(q)$ (respectively, $\GL(m,q^2) \times
G_{n-m}(q)$ if~$G$ is unitary). 

Let~$\St$ denote the character of the Steinberg representation of~$G$, and
let~$\omega$ be the class function introduced in Definition~\ref{DefinitionOmega}.
Then 
$$\om \cdot\,\St =\sum_{m=0}^n \left( \infl_{P_m} \left(\St^-_{m} 
\boxtimes \gamma'_{n-m} \right)\right)^G.$$
Here, $\St^-_{m} = 1^- \cdot\,\St_m$, where $\St_m$ denotes the Steinberg 
character of~$\GL(m,q)$ (respectively $\GL(m,q^2)$), and $1^-$ the unique 
linear character of this group of order~$2$, if ~$q$ is odd, and the trivial 
character, otherwise.

Moreover, $\gamma'_{n-m}$ is the Gelfand-Graev character of $G_{n-m}(q) = 
\GL(n-m,q)$ if $G$ is the general linear group. In the other cases, 
$\gamma'_{n-m}$ is a ``truncated'' Gelfand-Graev character of $G_{n-m}(q)$: It
is the sum of the regular characters of those Lusztig series which correspond
to semisimple elements without eigenvalue~$(-1)^q$ on~$V$.
\end{theor}
We are now going to discuss some consequences of the main result.
\begin{corol}\label{MultiplicityFree}
Let the notation be as in Theorem~\ref{Thm1Intro} and suppose that~$G$ is not
the general linear group. Then the character $\omega \cdot \St$ is multiplicity
free.
\end{corol}
We remark that this statement is not true for the general linear groups.

The above corollary is one of the principal ingredients in the proof of the 
following result. As indicated at the beginning of the introduction, this
also contains the main motivation for our work.

\begin{theor}\label{Thm2Intro}
Let~$G$ be one of the groups of Theorem~\ref{Thm1Intro} and let~$H'$ denote 
the stabilizer of a vector in the natural module for~$G$. Then the restriction
of the Steinberg character of~$G$ to~$H'$ is multiplicity free. In particular,
the same conclusion holds for the stabilizer~$H$ of a line.
\end{theor}
We do not know whether the analogous result holds for the orthogonal groups.

The irreducible characters of~$H'$ and~$H$ can be classified and our proof in
fact describes all the irreducible constituents of the restriction of~$\St$ to~$H'$
or~$H$ (see Subsection~\ref{RestrictingSteinberg}). In case~$G$ is a general 
linear group the above result is well known (see, e.g., \cite[Chapter~$5$]{BDK})
and its proof does not involve the product $\omega \cdot \St$.
To prove the result in case~$G$ is a unitary group and~$H'$ is the stabilizer
of an anisotropic vector (i.e., $H'$ is a unitary group of one degree less), 
we use in addition a nice observation by Olivier Brunat (see the appendix): 
The restriction of the Steinberg character of~$G$ to~$H'$ is the 
Weil-Steinberg character of~$H'$.

A result as in Theorem~\ref{Thm2Intro} is in general not true for other groups
of Lie type.  An example is provided by the Chevalley group $G_2(q)$. This 
group has two maximal standard parabolic subgroups~$P$ and~$Q$. Their character
tables have been computed in~\cite{AHu} in case $q$ is odd and not
a power of~$3$. Let~$q$ be such a prime power and let $G = G_2(q)$.
Then, in the notation of~\cite{AHu}, the restriction of $\St_G$ to~$P$
contains the irreducible character ${}_P\theta_2(0)$ with multiplicity
$(q + 1)/2$ (see \cite[Table A.4]{AHu}), and the restriction of $\St_G$ to~$Q$
has scalar product $q+1$ with the sum ${}_Q\theta_5(0) + {}_Q\theta_6(0)$ of
two irreducible characters (see \cite[Table A.7]{AHu}). So neither
is the restriction of $\St_G$ to the maximal parabolic subgroups multiplicity
free, nor are these multiplicities bounded independently of~$q$.

Theorem~\ref{Thm2Intro} has some interesting consequences for the $\ell$-modular 
representation theory of~$G$ for $\ell \nmid q$. Namely, the multiplicites of 
the $\ell$-modular constituents of (the reduction modulo~$\ell$) of the
Steinberg character of~$G$ can be controlled to some extent by the
$\ell$-modular decomposition numbers of~$P$. An example of such an
application to~$\Sp(6,q)$ is given in~\cite[Section~$5$]{AH}.

The Steinberg character is of defect~$0$ in the defining characteristic.
In this case, $\omega \cdot \St$ is the ordinary character of a projective
module~$M$. Thus Corollary~\ref{MultiplicityFree} yields PIMs which are 
multiplicity free as ordinary characters. In the case of the symplectic 
groups in characteristic~$2$ we were able to work out the decomposition of 
$M$ as a direct sum of PIMs.
In order to state this result, we need to recall some notions of
algebraic group theory. Let~$q$ be a power of~$2$ and let~$\mathbf{K}$ denote
an algebraic closure of the finite field~$\mathbb{F}_q$. Let $\mathbf{G} =
\Sp(2n,\mathbf{K})$ be the symplectic group of degree $2n$ over~$\mathbf{K}$.
Furthermore, let~$F$ be a standard Frobenius map of~$\mathbf{G}$, so that 
$G := \mathbf{G}^F = \Sp(2n,q)$ is the finite symplectic group of degree~$2n$
 over~$\mathbb{F}_q$ as in Theorem~\ref{Thm1Intro}.
If $\nu$ is a dominant weight of~$\mathbf{G}$ we denote by $\phi_\nu$ the 
rational irreducible representation of $\mathbf{G}$ corresponding to $\nu$. 
If $\nu$ is, furthermore, $q$-restricted, we write $\Phi_\nu$ for the principal
indecomposable character of~$G$ corresponding to the irreducible 
$\mathbb{F}_qG$-representation obtained by restricting $\phi_\nu$ to~$G$. 
\begin{theor} \label{Thm3Intro}
Let $\lambda_1 \ld \lambda_n$ be the fundamental weights of $\mathbf{G}$
(ordered as in Bourbaki \cite{Bo}).
Let $\nu_j = (q-1)\lambda_1+ \cdots +(q-1)\lambda_{n-1}+j\lambda_n$ for 
$0\leq j<q$. 
Then 
$$\omega \cdot \St = \sum_{j=0}^{q-1}\Phi_{\nu_j}.$$  
\end{theor}
It follows that the decomposition of every $\Phi_{\nu_j}$ as sum of ordinary
characters is multiplicity free.
We are not able to distribute the ordinary irreducible constituents of 
$\omega \cdot \St$ described in Theorem~\ref{Thm1Intro} between the 
projective indecomposable characters determined in Theorem \ref{Thm3Intro}. 
(This distribution will depend on the chosen $2$-modular system used to
define $\omega$ and the $\Phi_{\nu_j}$.)

Our approach is based on Deligne-Lusztig theory. In particular we have to pass
to dual groups in some arguments. The Weil characters of the classical groups
(where they exist) are closely related to properties of the natural module for
the groups. This is already apparent from the values of these characters on
semisimple elements as indicated above. Most important for our results, however, 
is the following property.  Consider a decomposition of the natural module into
a direct sum of non-degenerate subspaces. The stabilizer of this decomposition
is a direct product of classical groups induced on the subspaces, and the
Weil character restricts to this stabilizer as a product of the Weil
characters of these factors. Such stabilizers are in general not compatible
with duality of reductive groups. This is the reason why we take
some care in Sections~\ref{ToriInClassicalGroups} and~\ref{GeometricConjugacy}
to derive the necessary facts about maximal tori in duality and
their actions on the natural modules.

We conclude this introduction with an outline of the paper. In Section~$2$ 
we discuss maximal tori in classical groups and a decomposition of the natural 
module with respect to a given maximal torus. In Section~$3$ we relate these
decompositions for classical groups in duality. Section~$4$ introduces the
Weil representations and their characters and derives their properties needed
later on. In Section~$5$ we prove Theorem~\ref{Thm1Intro} for the symplectic and
unitary groups, as well as Corollary~\ref{MultiplicityFree}. The proof of 
Theorem~\ref{Thm1Intro} for the general linear groups is given in Section~$6$. 
It is different to the proof for the other classical groups. Section~$7$ is 
devoted to the applications of our main result, Theorems~\ref{Thm2Intro} 
and~\ref{Thm3Intro}.

\section{Tori in classical groups}
\label{ToriInClassicalGroups}

Let $V$ be a finite-dimensional non-degenerate unitary, symplectic or
orthogonal space over the finite field ${\mathbb F}_q$ with~$q$ elements
if~$V$ is symplectic or orthogonal, and $q^2$ elements if~$V$ is unitary.
We further assume that $\dim V$ is odd if~$V$ is orthogonal. In the
latter case we let~$G$ be the group of isometries of determinant~$1$,
otherwise~$G$ is the group of all isometries of~$V$. Thus $G$ is one of
the groups~$U(V)$, $\Sp(V)$, or $\SO(V)$.

In Subsection~\ref{TDecomposition} below we describe a decomposition of~$V$
relative to a maximal torus~$T$ of~$G$ and some formal properties of this
decomposition needed later on.

The concept of a maximal torus is defined via the algebraic group
underlying~$G$. We also have to compare such decompositions of the natural
module for groups which are dual to each other in the sense of Deligne and
Lusztig, with respect to dual maximal tori. In
Subsection~\ref{ToriClassification} we therefore introduce maximal tori and
the corresponding decompositions of~$V$ from an algebraic group point of view.
This treatment will also give proofs for the statements in~\ref{TDecomposition}
and allows us to avoid addressing uniqueness questions which arise for small
values of~$q$.

\subsection{The $T$-decomposition of~$V$}\label{TDecomposition}
Let $T$ be a maximal torus in $G$.

We will call an orthogonal direct sum decomposition
\begin{equation}
\label{EqTDecomposition}
V = V_0 \oplus V_1 \oplus \cdots \oplus V_k \oplus V_{k+1}
\oplus \cdots \oplus V_{k+l},
\end{equation}
a {\em $T$-decomposition of} $V$, if it has the following
properties:
\subsubsection{} \label{Cond211} The subspaces~$V_i$ are non-degenerate
$T$-submodules
for $1 \leq i \leq k + l$, $V_{k+1}, \ldots , V_{k+l}$ are irreducible
and $V_{1}, \ldots , V_{k}$ are reducible
and each of these $V_i$ is the sum of two irreducible, totally
singular $T$-submodules of equal dimension.
Moreover, $V_0 = \{0\}$ in the unitary and
symplectic case; otherwise $V_0$ is a $1$-dimensional subspace
spanned by an anisotropic vector, and $T$ acts trivially on~$V_0$.

\subsubsection{} \label{Cond212} For $1 \leq i \leq k + l$, let $G_i$ be the
subgroup of $G$
fixing $V_i$ and acting as the identity on the orthogonal complement of $V_i$.
Then $G_i \cong U(V_i)$, $\SO(V_i)$, or $\Sp(V_i)$, respectively. Let~$H$ be the
subgroup of~$G$ generated by the~$G_i$. Then~$H$ stabilizes all subspaces
$V_1,\ldots ,V_{k+l}$ and we have $H = G_1\times \cdots\times G_{k+l}$.
Put $T_i = T \cap {G_i}$. We then require that~$T_i$ is a cyclic maximal
torus of~$G_i$ for all~$i$ and
$$T = T_1 \times \cdots \times T_{k+l}.$$

\subsubsection{} \label{Cond213} Let $\mu_i= \dim V_i$ in case $V$ is a unitary
space. Then
$\mu_i$ is even for $1 \leq i \leq k$, and odd, otherwise. In the other cases,
each $V_i$ for $i \geq 1$ has even dimension and we write $\dim V_i = 2\mu_i$.
For $1 \leq i
\leq k$ we have $|T_i| = q^{\mu_i} - 1$, and for $k + 1 \leq i \leq k + l$ we
have $|T_i| = q^{\mu_i} + 1$.

\medskip

\noindent We will show below that a $T$-decomposition of~$V$ always exists.
Of course, the three conditions above are not independent. Clearly, one can
always find a decomposition~(\ref{EqTDecomposition}) of~$V$
satisfying~\ref{Cond211}. Also,~\ref{Cond213} follows
from~\ref{Cond212}, and, in a generic situation,~\ref{Cond212} is
implied by~\ref{Cond211}. Consider, however, the case~$q=2$
and~$V$ symplectic of dimension~$4$. Then $G = \Sp (4,2)$. Let~$T$
be the maximal torus of order~$3$ which is the Coxeter torus of
the split Levi subgroup $\GL(2,2)$ of~$G$. Thus there is a
$T$-decomposition of~$V$ with $k = 1$ and $l = 0$. There also is a
decomposition of~$V$ into an orthogonal direct sum of two
non-degenerate $2$-dimensional irreducible $T$-submodules. This
decomposition does not satisfy~\ref{Cond212}.


\begin{lemma} \label{ab5}
If~$G$ is unitary or $q > 2$, every maximal torus of $G$ induces a unique 
$T$-decomposition (up to reordering) of~$V$. 

Otherwise, any $T$-decomposition 
refines the decomposition $V = V^T \oplus ({V^T})^\perp$, where $V^T := 
\{v \in V \mid tv = v \text{\ for all\ } t\in T\}$. More precisely, 
$V^T = V_0 \oplus V_1 \oplus \ldots \oplus V_{k'}$ for some $k' \leq k$. 
The decomposition $({V^T})^\perp = V_{k'+1} \oplus \ldots \oplus V_{k+l}$ is 
unique (up to reordering), whereas the~$V_i$ in the decomposition $V^T = 
V_0 \oplus V_1 \oplus \ldots \oplus V_{k'}$ are hyperbolic planes (and so
this decomposition is not unique).
\end{lemma}
\begin{prf}
The existence of a $T$-decomposition will be proved in 
Subsections~\ref{UnitaryGroups} and \ref{SymplecticGroups} below. Suppose first
that $|T_i| > 1$. 
As~$T_i$ acts non-trivially on~$V_i$ but trivially on~$V_j$ for $j \neq i$, it
follows that $V_i$ and $V_j$ are not isomorphic (as $\mathbb{F}_{q^2}T$-modules 
respectively $\mathbb{F}_qT$-modules). Hence $V_i$ is a homogeneous component 
of~$V$ provided it is irreducible. 
Otherwise $V_i = V_i' \oplus V_i''$ and  $V_i'$, $V_i''$ are dual $T$-modules. 
If they are isomorphic,~$V_i$ is again a homogeneous component, and if they 
are not, each of $V_i'$, $V_i''$ is a homogeneous component of~$V$. 

Suppose now that $|T_i| = 1$. This can only happen if~$G$ is symplectic or 
orthogonal and $q = 2$ Then $V_i$ is acted on by $T$ trivially, $V_i$ is a 
hyperbolic plane, and $V^T$ is the sum of the $V_i$ with $T_i = 1$. This 
proves the assertions. In particular, the uniqueness statements follow from 
these observations.
\end{prf}


\subsection{Classification of maximal tori in finite reductive groups}
\label{ToriClassification}

Let~$q$ be a power of the prime~$p$, and let $\mathbf{K}$ denote an
algebraic closure of~$\mathbb{F}_p$.
We start with a connected reductive algebraic group~$\mathbf{G}$ over
$\mathbf{K}$, defined over $\mathbb{F}_q$, and denote by~$F$ the
corresponding Frobenius morphism. Closed, connected, $F$-stable subgroups
of $\mathbf{G}$ will be denoted by boldface letters, and if~$\mathbf{H}$
is such a subgroup, we write $H := \mathbf{H}^F := \{ h \in \mathbf{H}
\mid F(h) = h \}$ for the finite group of $F$-fixed points of
$\mathbf{H}$. The pair $(\mathbf{G},F)$, or simply the group $G =
\mathbf{G}^F$, is called a finite reductive group or a finite group of
Lie type.

To describe the maximal tori of~$G$ up to $G$-conjugacy, we follow
\cite[Section~$3.3$]{C}. Thus we fix an $F$-stable maximal torus
$\mathbf{T}_0$ of~$\mathbf{G}$, and let $W :=
N_{\mathbf{G}}(\mathbf{T}_0)/\mathbf{T}_0$ denote the corresponding Weyl
group of~$\mathbf{G}$. (Notice that the results of \cite[Section~$3.3$]{C}
are formulated for a maximally split torus $\mathbf{T}_0$, but that this
assumption is not needed; see \cite[3.23]{DM}.)

For every $w \in W$ we denote by $\dot{w}$ an element of $N :=
N_{\mathbf{G}}(\mathbf{T}_0)$ mapping to $w$ under the natural
epimorphism. For $t \in \mathbf{T}_0$ and $w \in W$ we let
$${^w\!t} := \dot{w}t\dot{w}^{-1}.$$
Clearly, the element ${^w\!t}$ does not depend on the particular
choice of~$\dot{w}$.

The $G^F$-classes of maximal tori in~$G$ are in bijection with the
$F$-conjugacy classes of~$W$. These are the orbits on $W$ under the
$F$-twisted $W$-action, also called $F$-conjugation, $w \mapsto
vwF(v)^{-1}, v, w \in W$.

This bijection arises as follows. Let $w \in W$.
By the Lang-Steinberg theorem, there is $g \in \mathbf{G}$
with $g^{-1}F(g) = \dot{w}$. Then $\mathbf{T} := {^g\mathbf{T}}_0$ is
$F$-stable and $T = \mathbf{T}^F = {^g(\mathbf{T}}_0^{{w}F})$,
where
$$\mathbf{T}_0^{{w}F} := \{ t \in \mathbf{T}_0 \mid {^w\!F}(t) = t \}.$$
Let $h \in \mathbf{G}$ with
$h^{-1}F(h) \in N$. Then ${^h\mathbf{T}}_0$ is $F$-stable and
${^h\mathbf{T}}_0$ is conjugate to ${^g\mathbf{T}}_0$ in~$G$ if and only
if the image of $h^{-1}F(h)$ in~$W$ is $F$-conjugate to $w$ in $W$.
We write $\mathbf{T}_w$ for any $F$-stable maximal torus of
$\mathbf{G}$ which corresponds to the $F$-conjugacy class of $w \in W$
in the way described above, and we say that $\mathbf{T}_w$
arises from $\mathbf{T}_0$ by twisting with~$w$.

Let $\mathbf{T}$ be an $F$-stable maximal torus of $\mathbf{G}$. We
put $W(\mathbf{T}) := W_{\mathbf{G}}(\mathbf{T}) :=
N_{\mathbf{G}}(\mathbf{T})/\mathbf{T}$ (so that
$W = W(\mathbf{T}_0)$). Then $F$ acts on $W(\mathbf{T})$, and we have
$W(\mathbf{T})^F \cong N_{\mathbf{G}}(\mathbf{T})^F/\mathbf{T}^F$ for
the set of $F$-fixed points on $W(\mathbf{T})$ (see
\cite[Section~$1.17$]{C}). If $\mathbf{T} = \mathbf{T}_w$ for some
$w \in W$, then $W(\mathbf{T})^F \cong C_{W,F}(w)$, the $F$-centralizer
of~$w$ (see \cite[Proposition $3.3.6$]{C}).

Let us write $\mathcal{S}(\mathbf{G})$ for the set of pairs $(\mathbf{T},s)$,
where $\mathbf{T}$ runs through the $F$-stable maximal tori of~$\mathbf{G}$
and $s \in T$. We are interested in classifying $\mathcal{S}(\mathbf{G})$ up
to $G$-conjugacy. For this purpose let
\begin{equation}
\label{SetP}
\mathcal{P} := \{ (w, t) \mid w \in W, t \in \mathbf{T}_0^{wF} \}.
\end{equation}
As indicated above, an element $(w,t) \in \mathcal{P}$ determines a
$G$-conjugacy class of elements of $\mathcal{S}(\mathbf{G})$.
The Weyl group~$W$ acts on $\mathcal{P}$ by $v.(w,t) := (vwF(v)^{-1},
{^v\!t})$ for $v \in W, (w,t) \in \mathcal{P}$. Two elements
of $\mathcal{P}$ are in the same $W$-orbit if and only if they determine
the same $G$-conjugacy class in $\mathcal{S}(\mathbf{G})$.

We will now give the specific examples to be used later on.

\subsection{The unitary groups}\label{UnitaryGroups}

Let $\mathbf{V}$ denote a vector space over $\mathbf{K}$ of dimension~$d$,
and fix a basis $v_1, v_2, \ldots, v_{d}$ of~$\mathbf{V}$. We then identify
$\mathbf{G} := \GL(\mathbf{V})$ with the matrix group $\GL(d,\mathbf{K})$.
To obtain the finite unitary groups, we let $F : \mathbf{G} \rightarrow
\mathbf{G}$ be the Frobenius morphism defined by $F(a_{ij}) :=
{((a_{ij}^q)^{-1})}^t$ for $(a_{ij}) \in \mathbf{G}$. Then $G = \mathbf{G}^F
= U(d,q) \leq \GL(d, q^2)$ with respect to the Hermitian form
$\sum_{i=1}^d x_iy_i^q$ on the $\mathbb{F}_{q^2}$-vector
space $V = \mathbf{V}(\mathbb{F}_{q^2})$ with basis $v_1, \ldots , v_d$.

In this case we choose $\mathbf{T}_0$ to be the group of diagonal matrices
of $\mathbf{G}$. (Thus $\mathbf{T}_0$ is not maximally split.)
Then $N = N_{\mathbf{G}}(\mathbf{T}_0)$ is the group of
monomial matrices and $W = N/\mathbf{T}_0$ can and will be identified
with the subgroup of permutation matrices of $\mathbf{G}$. Thus $W$ is
isomorphic to the symmetric group~$S_d$ on~$d$ letters, acting by permuting
the basis vectors $v_1 , \ldots , v_d$. Clearly,~$F$ acts trivially on~$W$.

The conjugacy classes of~$W$ are parametrized by the partitions of~$d$,
via the cycle type of a permutation. Let $w \in W$. Assume that~$w$
has~$k$ cycles of even lengths $\mu_1 \geq \cdots \geq \mu_k$, and~$l$
cycles of odd lengths $\mu_{k+1} \geq \cdots \geq \mu_{k+l}$.
We assume that $w = c_1 c_2 \cdots c_k c_{k+1} \cdots c_{k+l}$,
where $c_i$ is a cycle of length $\mu_i$. For $1 \leq i \leq k + l$, let
$\mathbf{V}_i$ denote the subspace of $\mathbf{V}$ spanned by the basis
vectors moved by $c_i$ (or by the unique basis vector corresponding to
$c_i$ if this is a $1$-cycle), and put $\mathbf{G}_i := \GL(\mathbf{V}_i)$.
The subspace $\mathbf{V}_i$ has
dimension~$\mu_i$ and $\mathbf{V} = \mathbf{V}_1 \oplus \cdots \oplus
\mathbf{V}_{k+l}$. We embed $\mathbf{G}_1 \times \cdots \times
\mathbf{G}_{k+l}$ into $\mathbf{G}$ in the natural way. Note that each
$\mathbf{G}_i$ is $F$-invariant, and that $G_i = \mathbf{G}_i^F \cong
U(\mu_i,q)$, acting on $V_i = \mathbf{V}_i(\mathbb{F}_{q^2})$, the
$\mathbb{F}_{q^2}$-subspace of $\mathbf{V}_i$ generated by
$\{v_1, \ldots, v_d \} \cap \mathbf{V}_i$.

Now choose $g_i \in \mathbf{G}_i$ with $g_i^{-1}F(g_i) = {c}_i$,
$1 \leq i \leq k+l$, and put $g := g_1 \times \cdots \times g_{k+l}$.
Then $g^{-1}F(g) = w$. Moreover, $\mathbf{T} := \mathbf{T}_w
:= {^g\mathbf{T}}_0 =
\mathbf{T}_1 \times \cdots \times \mathbf{T}_{k+l}$, with $\mathbf{T}_i
:= {^{g}\mathbf{T}}_{0,i} = {^{g_i}\mathbf{T}}_{0,i}$, where $\mathbf{T}_{0,i}
:= \mathbf{T}_0 \cap \mathbf{G}_i$, $1 \leq i \leq k+l$.
It follows that $T = \mathbf{T}^F = T_1 \times \cdots \times T_{s+t}$,
each $T_i$ acting on~$V_i$.

Fix $i$, $1 \leq i \leq k+l$, put $\mathbf{U} := \mathbf{V}_i$ and
$c:= c_i$. Let $u_1, \ldots , u_m$ be the basis vectors contained in
$\mathbf{U}$, numbered in such a way that $c$ maps $u_j$ to $u_{j+1}$,
$1\leq j \leq m$ (indices taken modulo~$m$). Write $h(\zeta_1, \ldots ,
\zeta_m)$ for the element of $\mathbf{T}_{0,i}$ which acts on $u_j$ by
multiplication with $\zeta_j \in \mathbf{K}^\times$, $1 \leq j \leq m$.
Then ${^c\!F}(h(\zeta_1, \ldots ,\zeta_m)) = h(\zeta_m^{-q}, \zeta_1^{-q},
\ldots, \zeta_{m-1}^{-q})$. Thus $h(\zeta_1, \ldots ,\zeta_m)$ is fixed
under the action of $cF$ if and only if $h(\zeta_1, \ldots ,\zeta_m) =
h(\zeta, \zeta^{-q}, \ldots , \zeta^{(-q)^{m-1}})$ for some $\zeta \in
\mathbf{K}$ with $\zeta^{(-q)^m} = \zeta$. It follows that $T_i$
is cyclic of order $q^m - 1$, if $m$ is even, and of order $q^m + 1$,
if $m$ is odd. In the former case,~$T_i$ fixes a maximal isotropic
subspace of~$V_i$, and in the latter case~$T_i$ acts irreducibly
on~$V_i$.

We have thus constructed a $T$-decomposition of~$V$ (see~\ref{TDecomposition}).

If $z \in C_{W,F}(w)$, then~$z$ permutes the cycles~$c_i$ of~$w$. Hence~$\dot{z}$
also permutes the tori $\mathbf{T}_{0,i}$, and so the corresponding
element $g \dot{z} g^{-1} \in W(\mathbf{T})^F$ permutes the tori~$\mathbf{T}_i$.

\subsection{The symplectic and orthogonal groups}\label{SymplecticGroups}
Let $\mathbf{V}$ be a vector space over $\mathbf{K}$ of dimension $d = 2n$
or $d = 2n + 1$. We choose a basis
\begin{equation}
\label{Basis}
v_1, v_2, \ldots, v_n, [v_0,] v_n', \ldots, v_2', v_1'
\end{equation}
of $\mathbf{V}$ (where $v_0$ is not present if $d = 2n$). The typical element
of $\mathbf{V}$ is denoted as $[x_0 v_0] + \sum_{i = 1}^nx_iv_i + x_i'v_i'$
with $x_0, x_i, x_i' \in {\mathbf{K}}$ (and without first summand if $d = 2n$).
Elements of $\GL(\mathbf{V})$ are written as matrices with respect to the
basis~(\ref{Basis}).

If $d = 2n$, we define a symplectic form on $\mathbf{V}$ such that
$v_i, v_{i}'$ is a hyperbolic pair for all $1 \leq i \leq n$ and such that
the planes $\langle v_i, v_i' \rangle$ are pairwise orthogonal. Let
$\mathbf{G} := \Sp( \mathbf{V} )$ denote the symplectic group with respect to
this form. We usually identify the elements of~$\mathbf{G}$ with their matrices
with respect to the basis~(\ref{Basis}), so that $\mathbf{G} = 
\Sp(2n, \mathbf{K}) \leq \GL(2n,\mathbf{K})$. We let~$F$ denote the
standard Frobenius morphism of $\mathbf{G}$ mapping the matrix
$(a_{ij})$ to $(a_{ij}^q)$. Then $G = \mathbf{G}^F = \Sp(2n,q) \leq
\GL(2n, q)$ with respect to the symplectic form $\sum_{i=1}^n (x_iy_{i}' -
x_{i}'y_i)$ on the $\mathbb{F}_{q}$-vector space $V = \mathbf{V}(\mathbb{F}_q)$
with basis $v_1, \ldots , v_{n}, v_n', \ldots, v_1'$.

If $d = 2n + 1$, we define the orthogonal form~$Q$ on 
$\mathbf{V}$ by $Q(x_0 v_0 + \sum_{i = 1}^nx_iv_i + x_i'v_i') := x_0^2 + 
\sum_{i=1}^nx_ix_i'$. Let $\mathbf{G} := \SO( \mathbf{V} ) = \SO(2n+1, 
\mathbf{K}) \leq \GL(2n+1,\mathbf{K})$ denote the special orthogonal group 
with respect to this form, and let $F$ be the
standard Frobenius morphism of $\mathbf{G}$. Then $G = \mathbf{G}^F =
\SO(2n+1,q) \leq \GL(2n+1, q)$ with respect to the orthogonal form
$x_0^2 + \sum_{i=1}^n x_ix_{i}'$ on the $\mathbb{F}_{q}$-vector
space $V = \mathbf{V}(\mathbb{F}_q)$ with basis $v_1, \ldots , v_{n},
v_0, v_n', \ldots, v_1'$.

Now let $\mathbf{V}$, $\mathbf{G}$, $F$ be one of the two configurations
introduced above. We choose $\mathbf{T}_0$ to be the group of diagonal
matrices of~$\mathbf{G}$. For $\zeta_1, \ldots , \zeta_n \in \mathbf{K}^\times$
we let $h(\zeta_1, \ldots , \zeta_n)$ denote the diagonal element of
$\mathbf{G}$ which acts by multiplication with $\zeta_i$ on $v_i$, and
by multiplication with $\zeta_i^{-1}$ on $v_{i}'$, $1 \leq i \leq n$.
Thus $\mathbf{T}_0 = \{ h(\zeta_1, \ldots , \zeta_n) \mid \zeta_1, \ldots ,
\zeta_n \in \mathbf{K}^\times \}$. (If $\mathbf{G}$ is orthogonal, every element
of $\mathbf{T}_0$ fixes~$v_0$.)

Let $W = N/\mathbf{T}_0$ with $N = N_{\mathbf{G}}(\mathbf{T}_0)$ denote
the Weyl group of $\mathbf{G}$. Then $W$ is the Weyl group of type $C_n$,
isomorphic to the wreath product of a cyclic group of order~$2$ with~$S_n$.
Clearly,~$F$ acts trivially on~$W$.

It is convenient to consider the faithful actions of~$W$ on the character
group $X := X(\mathbf{T}_0) := \Hom(\mathbf{T}_0, \mathbf{K}^\times)$ and on
the cocharacter group $Y := Y(\mathbf{T}_0) := \Hom(\mathbf{K}^\times,
\mathbf{T}_0)$ of~$\mathbf{T}_0$. These are free abelian groups of rank~$n$
with bases $\hat{e}_1, \ldots , \hat{e}_n$ defined by
$\hat{e}_i(h(\zeta_1, \ldots , \zeta_n)) = \zeta_i$, 
and $e_i(\zeta) = h(1, \ldots, 1, \zeta, 1, \ldots , 1)$
(where $\zeta$ is on position~$i$), respectively. The action of~$W$ on~$X$
and~$Y$ fixes the sets $\{ \pm \hat{e}_j \mid 1 \leq j \leq n \}$ and
$\{ \pm e_j \mid 1 \leq j \leq n \}$, respectively.

The set of conjugacy classes of $W$ is parametrized by the set of
bipartitions of~$n$. Let $w \in W$, viewed as a permutation group on
$\{ \pm \hat{e}_j \mid 1 \leq j \leq n \}$. Then~$w$ determines a bipartition
of~$n$ in the following way. There is a permutation $\pi = \pi(w)$ on
$\{1, \ldots ,n \}$ and a vector $(\sigma_1, \ldots , \sigma_n)$ of signs
(i.e., $\sigma_i \in \{+1, -1 \}$ for all $1 \leq i \leq n$) such that
${^w\hat{e}}_i = \sigma_i\hat{e}_{\pi(i)}$ for all $1 \leq i \leq n$.
The type of a cycle $(i_1, i_2, \ldots , i_m)$ of $\pi$ on $\{ 1, \ldots , n\}$
is the sign $\sigma_{i_1} \sigma_{i_2} \cdots \sigma_{i_m}$. Let
$\mu_1 \geq \mu_2 \geq \cdots \geq \mu_k$ denote the lengths of the
cycles of type~$+1$ of~$\pi$, and $\nu_1 \geq \nu_2 \geq \cdots \geq \nu_l$
the lengths of the cycles of type~$-1$ of~$\pi$. Then the pair $(\mu,\nu)$
with $\mu := (\mu_1, \ldots , \mu_k)$ and $\nu := (\nu_1, \ldots , \nu_l)$
is a bipartition of~$n$ which determines~$w$ up to conjugacy in~$W$.
Clearly, every bipartition of~$n$ arises in this way from a conjugacy
class of~$W$.

Let $w \in W$ correspond to the bipartition $(\mu,\nu)$ as above. Then
$w = c_1 c_2 \cdots c_k c_{k+1} \cdots c_{k+l}$ with pairwise commuting
elements $c_i \in W$, such that $\pi(c_i)$ is a cycle of type~$+1$ and
length $\mu_i$ for $1 \leq i \leq k$, and a cycle of type~$-1$ and length
$\nu_{i-k}$ for $k+1 \leq i \leq k + l$. The set of elements of
$\{ \pm \hat{e}_j \}$ moved by $c_i$ is invariant under multiplication
by~$-1$, and these sets form a partition of $\{ \pm \hat{e}_j \mid
1 \leq j \leq n \}$. We obtain a decomposition
$$X = X_1 \oplus \cdots \oplus X_{k+l}$$
into a direct sum of $w$-invariant, $w$-irreducible subgroups $X_i$
spanned by the orbits of $\langle w \rangle$ on $\{\pm \hat{e}_j \}$.
We have a corresponding decomposition
\begin{equation}
\label{YwDecomposition}
Y = Y_1 \oplus \cdots \oplus Y_{k+l}.
\end{equation}

For each $1 \leq i \leq k + l$, let $\mathbf{V}_i$ denote the subspace of
$\mathbf{V}$ spanned by the basis vectors corresponding to the
elements moved by $c_i$ (or to the two basis vectors $u, u'$ corresponding
to $c_i$ if this is a $1$-cycle), and put $\mathbf{G}_i := \Sp(\mathbf{V}_i)$
or $\mathbf{G}_i := \SO(\mathbf{V}_i)$, respectively.
The space $\mathbf{V}_i$ has dimension~$2\mu_i$ (with $\mu_i := \nu_{i-k}$
for $i > k$), and $\mathbf{V} = [\mathbf{V}_0 \oplus] \mathbf{V}_1 \oplus
\cdots \oplus \mathbf{V}_{k+l}$ (with $\mathbf{V}_0 := \langle v_0 \rangle$
in the orthogonal case). We embed $\mathbf{G}_1 \times \cdots \times
\mathbf{G}_{k+l}$ into $\mathbf{G}$ in the natural way. Note that each
$\mathbf{G}_i$ is $F$-invariant, and that $G_i = \mathbf{G}_i^F \cong
\Sp(2n_i,q)$ or $\SO^{\pm}(2n_i,q)$ (with $n_i = \mu_i$ or $\nu_i$),
acting on $V_i = \mathbf{V}_i(\mathbb{F}_q)$, the $\mathbb{F}_{q}$-subspace
of $\mathbf{V}_i$ generated by the basis vectors it contains.

Now choose $g_i \in \mathbf{G}_i$ with $g_i^{-1}F(g_i) = \dot{c}_i$,
$1 \leq i \leq k+l$, and put $g := g_1 \times \cdots \times g_{k+l}$.
If $1 \leq i \leq k$, the element ${c}_i$ lies in the stabilizer
of the maximal isotropic subspace generated by $v_1, \ldots , v_n$
and we choose $g_i \in \mathbf{G}_i$ also fixing this space.
Then $g^{-1}F(g) = w$. Moreover, $\mathbf{T} := {^g\mathbf{T}}_0 =
\mathbf{T}_1 \times \cdots \times \mathbf{T}_{k+l}$, with $\mathbf{T}_i
:= {^{g}\mathbf{T}}_{0,i} = {^{g_i}\mathbf{T}}_{0,i}$, where $\mathbf{T}_{0,i}
:= \mathbf{T}_0 \cap \mathbf{G}_i$, $1 \leq i \leq k+l$.
It follows that $T = \mathbf{T}^F = T_1 \times \cdots \times T_{k+l}$,
each $T_i$ acting on~$V_i$.

Fix $i$, $1 \leq i \leq k+l$, put $\mathbf{U} := \mathbf{V}_i$ and
$c:= c_i$. Let $\hat{e}_{j_1}, \ldots , \hat{e}_{j_m}$ be the elements
moved by $\pi(c)$, numbered in such a way that $\pi(c)$ maps
$\hat{e}_{j_r}$ to $\hat{e}_{{j_{r+1}}}$
$1\leq r \leq m$ (lower indices taken modulo~$m$). For $1 \leq r \leq m$,
put $u_r := v_{j_r}$ and $u_r' := v_{j_r}'$. Write $h(\zeta_1, \ldots ,
\zeta_m)$ for the element of $\mathbf{T}_{0,i}$ which acts on $u_j$ by
multiplication with $\zeta_j \in \mathbf{K}^\times$, $1 \leq j \leq m$.
Suppose first that $i \leq k$.
Then ${^c\!F}(h(\zeta_1, \ldots ,\zeta_m)) = h(\zeta_m^{q}, \zeta_1^{q},
\ldots, \zeta_{m-1}^{q})$. Thus $h(\zeta_1, \ldots ,\zeta_m)$
is fixed under the action of $cF$ if and only if $h(\zeta_1, \ldots ,\zeta_m) =
h(\zeta, \zeta^q, \ldots , \zeta^{q^{m-1}})$ for some $\zeta \in
\mathbf{K}^\times$ with $\zeta^{q^m} = \zeta$. It follows that $T_i$
is cyclic of order $q^m - 1$. Moreover,~$T_i$ fixes the maximal isotropic
subspace spanned by $v_1, \ldots , v_n$, by our choice of~$g_i$.
Next assume that $k+1 \leq i \leq k+l$. By conjugating $c = c_i$ by
a suitable element of $W$, we may and will assume that
${^c\!F}(h(\zeta_1, \ldots ,\zeta_m)) = h(\zeta_m^{-q}, \zeta_1^{q},
\ldots, \zeta_{m-1}^{q})$. Thus the $cF$-fixed points on $\mathbf{T}_0$
are of the form $h(\zeta, \zeta^q, \ldots , \zeta^{q^{m-1}})$ for some
$\zeta \in \mathbf{K}^\times$ with $\zeta^{q^m} = \zeta^{-1}$. Hence $T_i$ is
cyclic of order $q^m + 1$. Moreover,~$T_i$ acts irreducibly on~$V_i$.

Again, we have constructed a $T$-decomposition of~$V$.
As in the case of the unitary groups, we notice
that the elements of $W(\mathbf{T})^F$ permute the tori~$\mathbf{T}_i$.

\subsection{Neutral maximal tori}

We let $\mathbf{V}$, $\mathbf{G}$, $F$ be one of the configurations
introduced in~{\rm \ref{UnitaryGroups}} or~{\rm \ref{SymplecticGroups}},
and put $n=[(\dim \mathbf{V})/2]$ (the integer part). Thus $d = 2n$ or
$2n + 1$ in the situation of Subsection~\ref{UnitaryGroups} (and~$n$ has
the same meaning as in Subsection~\ref{SymplecticGroups} if $\mathbf{G}$
is symplectic or orthogonal).

We call a maximal torus $T = \mathbf{T}^F$ of~$G$ {\em neutral}, if no~$V_i$
in the $T$-decomposition of $V$, as specified above, is an irreducible
$T$-mod\-ule.

\begin{lemma}\label{n1}
{\rm (1)} If $(\mathbf{G},F)$ is as in~{\rm \ref{UnitaryGroups}} and if
$d = 2n + 1$ is odd, then ${G}$ does not have any neutral maximal torus.

{\rm (2)} Let $(\mathbf{G},F)$ be an orthogonal group as
in~{\rm \ref{SymplecticGroups}}, and let $\mathbf{T}$ be an $F$-stable
maximal torus of~$\mathbf{G}$. Consider the corresponding $T$-decomposition
$V = V_0 \oplus V_1 \oplus \cdots \oplus V_{k+l}$ of~$V$. If the Witt index
of $V_i$ is less than $(\dim V_i)/2$ for some $1 \leq i \leq k + l$, then
$i > k$. In particular,~$T$ is not neutral.

{\rm (3)}
If $(\mathbf{G},F)$ is as in~{\rm \ref{UnitaryGroups}}
or~{\rm \ref{SymplecticGroups}}, and if $d = 2n$ is even,
then the $G$-conjugacy classes of neutral maximal tori of~${G}$
are in a bijective correspondence with the set of partitions of~$n$.
\end{lemma}
\begin{prf}
We first prove (1) and (3). Suppose that we are in the situation
of~\ref{UnitaryGroups}. Then the torus $\mathbf{T}_w$ is neutral if and
only if $l = 0$, i.e., if and only if all $\mu_i$ are even.
In this case, $(\mu_1/2, \ldots , \mu_k/2)$ is a partition of $n = d/2$.

A torus $\mathbf{T}_w$ in the situation of~\ref{SymplecticGroups} is neutral
if and only if the partition $\nu$ is empty. Hence such tori are in bijection
with the set of bipartitions of $n$ of the form $(\mu, \mbox{\rm-} )$, where
$\mu$ runs through the partitions of~$n$.

To prove~(2), observe that for $i \leq k$, an irreducible $T$-submodule of
$V_i$ is maximal singular of dimension $(\dim V_i)/2$.
\end{prf}

\begin{lemma}\label{n2}
Let $(\mathbf{G},F)$ be as in~{\rm \ref{UnitaryGroups}}
or~{\rm \ref{SymplecticGroups}} and let~$\mathbf{T}$ be a neutral maximal
torus in $\mathbf{G}$ corresponding to the partition $(1^{m_1}, 2^{m_2},
\ldots ,n^{m_n})$ of~$n$. Then
$|W(\mathbf{T})^F|= \prod_{i=1}^n (2i)^{m_i} m_i!$.
\end{lemma}
\begin{prf}
This follows directly from $W(\mathbf{T}_w)^F \cong C_{W,F}(w)$
and the well known descriptions of the $F$-centralizers in the
respective Weyl groups.
\end{prf}

\subsection{Some notation}\label{SomeNotation}
We end this section by introducing some character theoretic notation, where the
word character refers to a complex character of a finite group. Let $X$ and $Y$
be finite groups. We denote by $\rho_X$ and $1_X$ the regular and the trivial
character of $X$. If $X$ has a unique cyclic quotient group of even order, we
denote by $1^-_X$ the non-trivial linear character of $X$ with values $\pm 1$.
For uniformity of some expressions, if $X$ is of odd order, we interpret
$1^-_X$ as $1_X$. If $\chi$ and $\psi$ are characters of $X$ and $Y$,
respectively, $\chi \boxtimes \psi$ denotes their outer product, a
character of $X \times Y$. In contrast, we use the symbol $\otimes$ to denote
the (inner) tensor product of representations of~$X$. If $Y$ is a subgroup 
of~$X$, then $\chi_Y$ is the
restriction of~$\chi$ to~$Y$, and $\psi^X$ the character of~$X$ induced from
$\psi$. Finally the usual inner product of two complex class functions $\chi$
and $\psi$ of~$X$ is denoted by $(\chi,\psi)$.

\section{Duality and geometric conjugacy}\label{GeometricConjugacy}

Let $(\mathbf{G}, F)$ be a finite reductive group.
We have to investigate
the dual reductive group $(\mathbf{G}^*, F^*)$ to some extent. In particular,
we wish to describe the pairs $(\mathbf{T}, \theta)$, where $\mathbf{T}$ is a
maximal $F$-stable torus of $\mathbf{G}$, and $\theta$ is an irreducible
(complex) character of $T$, up to conjugation in~$G$. This is most
conveniently done by passing to the dual group.
We fix a maximal $F$-stable torus $\mathbf{T}_0$ of $\mathbf{G}$, and a
maximal $F^*$-stable torus $\mathbf{T}_0^*$ of $\mathbf{G}^*$ satisfying the
conditions of \cite[Proposition 4.3.1]{C}.  In other words, $(\mathbf{G}, F)$
and $(\mathbf{G}^*, F^*)$ are in duality with respect to the pair
$(\mathbf{T}_0, \mathbf{T}_0^*)$. Again, the assumption of~\cite{C},
that the tori be maximally split, is not needed.
In the following, we mark the objects associated with $\mathbf{G}^*$ with
an asterisk.

\subsection{Geometric conjugacy.}\label{Section31}
We identify
$X := \Hom( \mathbf{T}_0, \mathbf{K}^\times )$ with
$Y^* := \Hom( \mathbf{K}^\times, \mathbf{T}_0^* )$ and
$Y := \Hom( \mathbf{K}^\times, \mathbf{T}_0 )$ with
$X^* := \Hom( \mathbf{T}_0^*, \mathbf{K}^\times )$.
Denote by $W := W(\mathbf{T}_0)$ and $W^* := W(\mathbf{T}^*_0)$ the
Weyl groups of~$\mathbf{G}$ and of~$\mathbf{G}^*$, respectively.

The identification of~$Y$ with~$X^*$ yields an $F$-$F^*$-equivariant
isomorphism
$$\delta: \Hom(Y,\mathbf{K}^\times) = \Hom(X^*,\mathbf{K}^\times)
\rightarrow \mathbf{T}_0^*$$
of abelian groups. For the isomorphism $\Hom(X^*,\mathbf{K}^\times)
\rightarrow \mathbf{T}_0^*$ see \cite[Propostion~3.1.2(i)]{C}.
As in
\cite[Proposition~4.2.3]{C}, there is an anti-isomorphism $W \rightarrow W^*,
w \mapsto w^*$, such that $\delta({^{w^{-1}}\!\psi}) = {^{w^*}\!\delta}(\psi)$
for all $\psi \in \Hom(Y,\mathbf{K}^\times)$ and $w \in W$.

Put
$$\mathcal{Q} := \{ (w, \psi) \mid w \in W, \psi \in \Hom(Y,\mathbf{K}^\times),
F({^{w^{-1}}\!\psi}) = \psi \}.$$
Then $W$ acts on $\mathcal{Q}$ by $v.(w,\psi) := (vwF(v)^{-1}, {^v\!\psi} )$
for $v \in W, (w,\psi) \in \mathcal{Q}$,  and there is a bijection
\begin{equation}
\label{QPBijection}
\mathcal{Q} \rightarrow \mathcal{P}^*,\quad\quad (w, \psi) \rightarrow
(F^*(w^*), \delta(\psi)).
\end{equation}
(For the definition of $\mathcal{P}^*$ see~(\ref{SetP}).)
One easily checks that $v.(w, \psi)$ is mapped to
${v^*}^{-1}.(F^*(w^*), \delta(\psi)).$
In particular, this map induces a bijection of the $W$-orbits in $\mathcal{Q}$
with the $W^*$-orbits in $\mathcal{P}^*$.

Let us write $\mathcal{T}(\mathbf{G})$ for the set of pairs
$(\mathbf{T},\theta)$, where $\mathbf{T}$ runs through the $F$-stable
maximal tori of~$\mathbf{G}$ and $\theta \in \Irr(T)$.

An element of
$\mathcal{Q}$ gives rise to a $G$-conjugacy class of elements of
$\mathcal{T}(\mathbf{G})$ as follows.
Choose an isomorphism
$$\Omega_{p'} \rightarrow \mathbf{K}^\times,$$
where $\Omega_{p'} \subseteq \mathbb{C}$ denotes the set of roots of unity
of $p'$-order (see \cite[Proposition~$3.1.3$]{C}).
Let $(w, \psi) \in \mathcal{Q}$. The condition $F({^{w^{-1}}\!\psi}) = \psi$ is
equivalent to $(wF - \mbox{\rm id})Y \leq \ker(\psi)$. Hence
$\psi$ may be viewed as an element of $\Hom( Y/(wF - \mbox{\rm id})Y,
\mathbf{K}^\times) \cong \Hom( Y/(wF - \mbox{\rm id})Y, \Omega_{p'})$.
Moreover, $Y/(wF - \mbox{\rm id})Y \cong \mathbf{T}_0^{wF}$ (see
\cite[Proposition 3.2.2]{C}). We thus obtain a pair $(\mathbf{T},\theta)
\in \mathcal{T}(\mathbf{G})$ with $\mathbf{T} = \mathbf{T}_w$ and where
$\psi$ is related to $\theta$ via an isomorphism
\begin{equation}
\label{PsiAndTheta}
\Hom( Y/(wF - \mbox{\rm id})Y, \mathbf{K}^\times) \rightarrow
\Hom( \mathbf{T}_0^{wF}, \Omega_{p'}) =
\Hom( \mathbf{T}_0^{wF}, \mathbb{C}^\times ).
\end{equation}
This construction yields a one-to-one correspondence between the set of
$W$-orbits on $\mathcal{Q}$ and the set $G\backslash\mathcal{T}(\mathbf{G})$
of $G$-conjugacy classes on
$\mathcal{T}(\mathbf{G})$. Through the bijection~(\ref{QPBijection})
and the considerations in~\ref{ToriClassification}, we obtain
a one-to-one correspondence
\begin{equation}
\label{TGandSG}
G\backslash\mathcal{T}(\mathbf{G}) \rightarrow
G^*\backslash\mathcal{S}(\mathbf{G}^*),
\end{equation}
where $G^*\backslash\mathcal{S}(\mathbf{G}^*)$ denotes the set of
$G^*$-conjugacy classes on $\mathcal{S}(\mathbf{G}^*)$.
We say that $(\mathbf{T}, \theta) \in \mathcal{T}(\mathbf{G})$ and
$(\mathbf{T}^*,s^*) \in \mathcal{S}(\mathbf{G}^*)$ are {\em dual}, if their
respective conjugacy classes correspond via~(\ref{TGandSG}).

Finally, the bijection~(\ref{QPBijection}) yields an isomorphism
$$\Irr( \mathbf{T}_w^F ) \rightarrow {\mathbf{T}^*_{F^*(w^*)}}^{F^*}$$
for every $w \in W$.

For $(\mathbf{T},\theta) \in \mathcal{T}(\mathbf{G})$ we put
$W(\mathbf{T})_\theta^F := \{ w \in W(\mathbf{T})^F \mid {^w\!\theta} =
\theta \}$ (for the definition of $W(\mathbf{T})$ see
Subsection~\ref{ToriClassification}). Similarly, if $(\mathbf{T}^*,s^*)
\in \mathcal{S}(\mathbf{G}^*)$, we put $W(\mathbf{T}^*)_{s^*}^{F^*}
:= \{ w \in W(\mathbf{T}^*)^{F^*} \mid {^w\!s}^* = s^* \}$.

We will need the following lemma later on.
\begin{lemma}
\label{WeylGroupFactoriazation}
Let $(\mathbf{G},F)$ be a unitary group as in
Subsection~{\rm \ref{UnitaryGroups}} or a symplectic group as in
Subsection~{\rm \ref{SymplecticGroups}}, and let $(\mathbf{T},\theta)
\in \mathcal{T}(\mathbf{G})$. Consider a $T$-decomposition of $V$ as
constructed in these subsections.

Put $I := \{ 1 \leq i \leq k + l \mid \theta_i = 1^-_{T_i} \}$, and
$J := \{ 1 \leq i \leq k + l \mid \theta_i \neq 1^-_{T_i} \}$.
Next, let $\mathbf{V}_I := \oplus_{i \in I} \mathbf{V}_i$,
and $\mathbf{V}_J := \oplus_{i \in J} \mathbf{V}_i$, so that
$\mathbf{V} = \mathbf{V}_I \oplus \mathbf{V}_J$.

Then the stabilizer in $\mathbf{G}$ of this decomposition equals
$\mathbf{G}_I \times \mathbf{G}_J$, where $\mathbf{G}_I$ and
$\mathbf{G}_J$ act as the identity on $\mathbf{V}_J$ and
$\mathbf{V}_I$, respectively. Moreover $\mathbf{T} = \mathbf{T}_I
\times \mathbf{T}_J$ with the $F$-stable tori $\mathbf{T}_I :=
\mathbf{T} \cap \mathbf{G}_I$ and  $\mathbf{T}_J := \mathbf{T}
\cap \mathbf{G}_J$.

Put $\theta_I := \theta_{T_I}$ and $\theta_J := \theta_{T_J}$.
Then
$$W(\mathbf{T})^F_\theta = W_{\mathbf{G}_I}(\mathbf{T}_I)^F_{\theta_I}
\times W_{\mathbf{G}_J}(\mathbf{T}_J)^F_{\theta_J}.$$
\end{lemma}
\begin{prf}
First note that the stabilizer of the orthogonal decomposition
$\mathbf{V} = \mathbf{V}_I \oplus \mathbf{V}_J$ equals
$\mathbf{G}_I \times \mathbf{G}_J$, since $\mathbf{G}$ is a
general linear or a symplectic group. Let $w \in
W(\mathbf{T})^F_\theta$, and choose an inverse image $\dot{w} \in
N_G(\mathbf{T})$ of~$w$. Since~$\dot{w}$ fixes~$\theta$, and
since~$\dot{w}$ permutes the factors ${T}_i$ of~${T}$ by the final
remarks of Subsections~\ref{UnitaryGroups}
and~\ref{SymplecticGroups}, it follows that~$\dot{w}$ normalizes
$T_I$ and $T_J$.

Now $|T_j| > 1$ for each $j \in J$ 
and if $T_j = \langle t_j \rangle$, then $t_j$ does not have eigenvalue~$1$
on~$V_j$. This implies that $V_I := \sum_{i\in I}V_i$ equals the fixed space
of~$T_J$. 

Since~$\dot{w}$ normalizes $T_J$, it follows 
that~$\dot{w}$ fixes $V_I$ and thus also $V_J = V_I^\perp$, and in turn it
fixes $\mathbf{V}_I$ and $\mathbf{V}_J$. Thus $\dot{w}$ is contained in
$\mathbf{G}_I \times \mathbf{G}_J$.

Hence $\dot{w} \in (\mathbf{G}_I \times \mathbf{G}_J)^F = \mathbf{G}_I^F
\times \mathbf{G}_J^F$, and so $\dot{w} = \dot{w}_I \cdot \dot{w}_J$
with $\dot{w}_I \in N_{\mathbf{G}_I}(\mathbf{T}_I)^F$ and
$\dot{w}_J \in N_{\mathbf{G}_J}(\mathbf{T}_J)^F$. Writing $w_I$ and~$w_J$
for the images of $\dot{w}_I$ and $\dot{w}_J$ in $W_\mathbf{G}(\mathbf{T})^F$,
respectively, we obtain $w_I \in W_{\mathbf{G}_I}(\mathbf{T}_I)^F_{\theta_I}$
and $w_J \in W_{\mathbf{G}_J}(\mathbf{T}_J)^F_{\theta_J}$, and hence the result.
\end{prf}

\subsection{Duality and $T$-decompositions.}\label{Section32}

Let $(\mathbf{G},F)$ be a unitary group as in~\ref{UnitaryGroups},
or a symplectic group as in~\ref{SymplecticGroups}.
If $(\mathbf{G},F)$ is the finite unitary group as
in~\ref{UnitaryGroups}, we may and will identify $(\mathbf{G}, F)$ with
its dual $(\mathbf{G}^*, F^*)$ and put $\mathbf{T}_0 = \mathbf{T}_0^*$.
If $(\mathbf{G},F)$ is the symplectic group as in~\ref{SymplecticGroups},
then $(\mathbf{G}^*,F^*)$ is the special orthogonal group of dimension
$2n+1$, also described in~\ref{SymplecticGroups}. As our reference torus
$\mathbf{T}_0^*$ in $\mathbf{G}^*$ we take the torus denoted by $\mathbf{T}_0$
in~\ref{SymplecticGroups}.

\begin{lemma}
\label{SameOrders}
Let~$(\mathbf{G},F)$ be a unitary group as in~{\rm \ref{UnitaryGroups}}
or a symplectic group as in~{\rm \ref{SymplecticGroups}}.
Suppose that $(\mathbf{T}, \theta) \in \mathcal{T}(\mathbf{G})$ and
$(\mathbf{T}^*,s^*) \in \mathcal{S}(\mathbf{G}^*)$ are dual pairs.
Then the following statements hold.

{\rm (a)}
$W(\mathbf{T})^F_\theta \cong W(\mathbf{T}^*)^{F^*}_{s^*}$.

{\rm (b)} Let
$$V = V_1 \oplus \cdots \oplus V_k \oplus V_{k+1} \oplus \cdots \oplus V_{k+l}$$
be a $T$-decomposition of~$V$ as constructed in~{\rm \ref{UnitaryGroups}}
or~{\rm \ref{SymplecticGroups}}.
Then there is a corresponding $T^*$-decomposition
$$V^* = V^*_0 \oplus V^*_1 \oplus \cdots \oplus V^*_k \oplus V^*_{k+1} \oplus
\cdots \oplus V^*_{k+l},$$
of~$V^*$ with $\dim V_i = \dim V_i^*$ for $1 \leq i \leq k + l$.

Consider the induced direct decompositions
$$(\mathbf{T}_1 \times \cdots \times \mathbf{T}_{k+l},\theta_1 \boxtimes
\cdots \boxtimes \theta_{k+l})$$
of $(\mathbf{T},\theta)$ and
$$(\mathbf{T}^*,s^*) = (\mathbf{T}^*_1 \times \cdots \times
\mathbf{T}^*_{k+l}, s_1^* \times \cdots \times s_{k+l}^*)$$
of $(\mathbf{T}^*,s^*)$.
Then the order of $\theta_i \in \Irr(T_i)$ equals the order of $s_i^*$
as automorphism on~$V^*_i$, for $1 \leq i \leq k + l$.
In particular, $\theta_i = 1^-_{T_i}$ if and only if $s^*_i$ acts as
$-1$ on $V^*_i$. Similarly, $\theta_i = 1_{T_i}$ if and only if $s^*_i$
acts as the identity on $V^*_i$.
\end{lemma}
\begin{prf}
The isomorphism in~(a) is derived in \cite[p.~289]{C}.

By conjugating in $G$ and $G^*$, respectively, we may assume that
$(\mathbf{T},\theta)$ is constructed from $(w, \psi) \in
\mathcal{Q}$ as in Subsection~\ref{Section31} and that
$(\mathbf{T}^*,s^*)$ corresponds to $(w^*, \delta(\psi))$ as in
Subsection~{\rm \ref{ToriClassification}}. (We remark that~(a) now
also follows from the fact that $W(\mathbf{T})^F_\theta$ and
$W(\mathbf{T}^*)^{F^*}_{s^*}$ are isomorphic to the stabilizers of
the pairs $(w,\psi) \in \mathcal{Q}$ and $(w^*,\delta(\psi))$,
respectively.) Notice that the conjugacy classes of $w$ and of
$w^*$ are labelled by the same partition, respectively bipartition
(since inverse elements are conjugate). We construct $\mathbf{T} =
\mathbf{T}_w$, $\mathbf{T}^* = \mathbf{T}_{w^*}^*$ and the
corresponding decompositions of~$V$ and~$V^*$ as
in~\ref{UnitaryGroups} and~\ref{SymplecticGroups}, respectively.
Considering the decompositions~(\ref{YwDecomposition}) of~$Y$
arising from~$w$, and of~$\mathbf{T}^*_0$ arising from $w^*$, we
obtain the following commutative diagram of abelian groups. $$
\xymatrix{ \Hom( Y, \mathbf{K}^\times ) \ar[d]^\delta
\ar[r]^-{\beta} & \Hom( Y_1, \mathbf{K}^\times ) \oplus \cdots
\oplus \Hom( Y_{k+l}, \mathbf{K}^\times ) \ar[d]^{\times \delta_i}
\\ \mathbf{T}_0^* \ar[r]^-{\beta^*} & \mathbf{T}_{0,1}^* \times
\cdots \times \mathbf{T}_{0,k+l}^* \\ } $$ All isomorphisms are
compatible with the actions of $\langle w, F \rangle$ in the top
row and $\langle w^*, F^* \rangle$ in the bottom row. Writing
$\beta( \psi ) = \sum_{i=1}^{k+l} \psi_i$ with $\psi_i \in \Hom(
Y_i, \mathbf{K}^\times )$, the characters $\theta_i$ correspond to
$\psi_i$ and the elements $s_i^*$ correspond to $\delta_i(\psi_i)$
under the group isomorphisms~(\ref{PsiAndTheta}). This gives the
first result.

Finally, as $V_i^*$ has no proper non-degenerate $T_i^*$-invariant subspace,
$s^*_i$ has order~$2$ if and only if it acts as~$-1$ on~$V^*_i$.
Since~$T_i$ is cyclic, the element $\theta_i \in \Irr(T_i)$ has order~$2$
if and only if $\theta_i = 1^-_{T_i}$. The last statement is trivial.
This completes the proof.
\end{prf}

\section{The characters of the Weil representations}
\label{WeilCharacters}

\subsection{The ordinary case}\label{OrdinaryCase}
Let $G = \Sp(2n,q)$ with $q$ odd, or $U(d,q)$, with $q$ arbitrary. Let
$V$ denote the natural module for $G$ and let $T$ be a maximal torus
of $G$.

The standard reference for Weil representations is G{\'e}rardin~\cite{G},
who computed their characters.
If $G = U(d,q)$, there is a unique  Weil representation of~$G$ (up to
equivalence). If $G=\Sp(2n,q)$, there are two  Weil representations of~$G$
(see \cite[Theorem 2.4(d)]{G}),
but the character values of the two Weil representations
on semisimple elements are the same (see \cite[Corollary 4.8.1]{G}).

Let $\hat{\omega}^{(G)}$ denote the character of a Weil representation
of~$G$. If~$G$ is symplectic, we put $\omega := \omega^{(G)} :=
\hat{\omega}^{(G)}$, and if~$G$ is unitary, we put $\omega := \omega^{(G)}
:= 1_G^- \cdot \hat{\omega}^{(G)}$. (Thus in the latter case,~$\omega$
is not the character of G{\'e}rardin's Weil representation if~$q$ is
odd.)

The most important feature of the Weil representation is the
multiplicative nature of its character. Namely, if $V = U \oplus U'$ where
$U$ and $U'$ are non-degenerate and mutually orthogonal then the embedding
$H := G_{U} \times G_{U'}$ into~$G$ gives $\om_{H} = \om^{(G_U)} \boxtimes
\om^{(G_{U'})}$ (see \cite[Corollaries 2.5, 3.4]{G}).

Let~$T$ be a maximal torus of~$G$. Corresponding to a
$T$-decomposition of~$V$ we have an induced
decomposition $T = T_1 \times \cdots \times T_{k+l}$ of~$T$, and a
subgroup $H = G_1\times \cdots \times G_{k+l}$ of~$G$. The above implies
that $\omega_H = \omega_{1} \boxtimes \cdots \boxtimes \omega_{{k+l}}$,
with $\omega_i := \omega^{(G_i)}$, $1 \leq i \leq k+l$.

\begin{lemma}\label{g1}
{\rm (a)} Suppose that $k=1$, $l=0$. Then $\omega_T = \rho_T+1^-_T$.

{\rm (b)} Suppose that $k=0$, $l=1$. Then $\omega_T=\rho_T-1^-_T$.

{\rm (c)} In general, we have
$$\om_T =
(\rho_{T_1}+1^-_{T_1})\boxtimes\cdots
\boxtimes(\rho_{T_k}+1^-_{T_k})\boxtimes (\rho_{T_{k+1}}
-1^-_{T_{k+1}})\boxtimes\cdots \boxtimes(\rho_{T_{k+l}}-1^-_{T_{k+l}}) .$$
\end{lemma}
\begin{prf}
The statements in (a) and (b) can be derived from
\cite[Corollaries 4.8.1, 4.8.2]{G}. The last statement follows from these.
\end{prf}

\subsection{The modular case}\label{ModularCase}
We change the point of view and consider instead the $p$-modular version of 
the Weil representation. If $p>2$ this is just the Brauer reduction modulo~$p$
of the Weil representation. If $G = \Sp(2n,q)$ and $p=2$ the Weil representation
does not exist, but it has been shown by the second author in \cite{Z90}, that
the analogue of its Brauer reduction modulo~$2$ does exist, and that this is 
exactly the generalized spinor representation of~$G$. If $q=2$, this is the 
usual spinor representation. 

Let~$q$ be a power of~$2$ and let~$\mathbf{K}$ denote an algebraic closure of
the finite field~$\mathbb{F}_q$. Let $\mathbf{G} = \Sp(2n,\mathbf{K})$ be the
symplectic group of degree $2n$ over~$\mathbf{K}$ as introduced in 
Subsection~\ref{SymplecticGroups} and let~$F$ be the standard Frobenius map 
of~$\mathbf{G}$ raising every matrix entry of~$\mathbf{G}$ to its $q$th power. 

To introduce the generalized spinor representation of $G = \Sp(2n,q)$, we 
recall some notions of algebraic group theory. 
Let $\lambda_1 \ld \lambda_n$ be the fundamental weights of $\mathbf{G}$ 
(ordered as in Bourbaki \cite{Bo}). An integer linear combination $\sum
a_i \lambda_i$ is called a weight of $\mathbf{G}$, and the weights
with $a_i \geq 0$ for $i=1\ld n$ are called dominant. There is a
canonical bijective correspondence between the dominant weights
and the equivalence classes of rational irreducible representations of
$\mathbf{G}$, and for a dominant weight $\nu$ we denote by
$\phi_\nu$ the irreducible representation of $\mathbf{G}$ corresponding 
to $\nu$. We set $\sigma_n = {\left(\phi_{(q-1)\lambda_n}\right)}_{G}$ and call
$\sigma_n$ the generalized spinor representation of $G$, while the spinor
representation is ${(\phi_{\lambda_n})}_{G}$. 
To avoid confusion we sometimes use the notation $\sigma_{n,q}$
for $\sigma_n$.

\subsubsection{The  Weil representation of the extrasymplectic group} Despite 
the fact that the representation $\sigma_n$ is explicitly constructed, its
Brauer character does not seem to have been computed. We need to do this
and, moreover, to express it in terms of characters of the maximal
tori in $G$. We could do this by straightforward computations
but it is more conceptual  to connect this with complex
representations of extraspecial $2$-groups.

\def\ex{extraspecial }
So we start with extraspecial groups. For a natural number $n$
there are two \ex groups of order $2^{2n+1}$ which we denote by
$E_n^+$ and $E_n^-$. The center $Z$ of each of them is of order 2.
  The central quotients are elementary
abelian 2-groups.  Let $C_4$ denote the cyclic group of order 4 and
let $E_n$ be the central product $C_4\cdot E_n^+$ (with common
subgroup of order 2). Then $C_4\cdot E_n^+ =C_4\cdot E_n^-$, so
$E_n$ contains $E_n^+$ and $E_n^-$ as subgroups of index 2. We
denote the central quotient by $V_n$ in all three cases. Then
the mapping $xZ \mapsto x^2$ defines a non-degenerate quadratic form on
$V_n$ and the two forms corresponding to $E_n^+$ and $E_n^-$ are
non-equivalent. The mapping $xZ \times yZ \mapsto [x,y]$ for $x,y\in
E_n$ defines a non-degenerate alternating form on $V_n$ which is
the polarization of both quadratic forms. Details   can be
found in \cite[page 80]{DH}.  Furthermore, $\Aut E_n^+/\Inn
E_n^+\cong O^+(2n,2)$,  $\Aut E_n^-/\Inn E_n^-\cong O^-(2n,2)$
(\cite[Theorem 20.8]{DH}) and $\Aut E_n/\Inn E_n\cong
\Sp(2n,2)\times C_2$. We denote by  $\Aut^0 E_n$ the subgroup of
$\Aut E_n$ consisting of the automorphisms acting trivially on the
center. So $\Aut^0 E_n/\Inn E_n \cong \Sp(2n,2)$.

It is also well known that every faithful complex irreducible representation of
$E_n$ has degree $2^n$, and its character $\chi$ vanishes on all non-central
elements. As elements of the center of $E_n$ are represented by scalar 
matrices, there are exactly two non-equivalent faithful irreducible 
representations of $E_n$ which are dual to each other. We denote any one of 
them by $\eta$. Let $\al$ be an automorphism of $E_n$ acting trivially on the
center. Then $\eta^\al = \eta$. It follows that $\eta (\al (x)) = g \eta(x) 
g^{-1}$ for some $g \in \GL(2^n,\mathbb{C})$. As~$g$ is determined by $\al$ up
to a scalar multiple, the mapping $\Aut^0 E_n\ra \GL(2^n,\mathbb{C})$ obtained
from this provides a projective representation $\pi$ of $\Aut^0 E_n$ into
$\GL(2^n,\mathbb{C})$. An irreducible projective representation of a finite 
group can be obtained from an ordinary representation of a central extension. 
It turns out that a central extension of $\Aut^0 E_n$ by a cyclic group of
order~$4$ is sufficient. Thus, there 
exists a group $R = R(n,2)$ with normal subgroup $E_n$ such that $R/E_n\cong
\Sp(2n,2)$, and an irreducible representation $\eta$ of $R$ of degree $2^n$
such that $\eta_{E_n}$ is irreducible. 

It is well known that the group $\Sp(2m,2^k)$ is isomorphic to a subgroup of
$\Sp(2mk,2)$. We fix an embedding $\Sp(2m,2^k)\ra \Sp(2mk,2)$ and denote by
$\ESp(2m,q)$ for $q=2^k$ the preimage of $\Sp(2m,q)$ in $R=R(mk,2)$. We call
$\ESp(2m,q)$ the extrasymplectic group and use the term ``Weil character'' 
for the character of its irreducible representation of degree $2^n = q^m$. 
The Weil character depends on $\eta$ which is immaterial for what follows as
we are only interested in the values of  $\eta$ at odd order elements. These
are independent of the choice of $\eta$.

\medskip
\begin{remar}
{\rm
(1) Usually the Weil character is considered for symplectic groups in odd
characteristic. However, there is a strong similarity between the odd 
characteristic Weil character at semisimple elements and the above introduced
Weil character for the extrasymplectic group at semisimple elements. Observe
that $\ESp(2m,q)$ is not split over $E_n$ so one cannot restrict $\eta$ to 
$\Sp(2m,q)$ in contrast to the case of odd $q$.

(2) The existence of the above projective representation of $\Aut^0 E_n$ was
probably shown first in Suprunenko \cite[Theorem 11]{Su} but he deals with the
linear group $\eta (E_n)\cdot S$ where $S$ is the group of all non-zero scalar
matrices. The observation that the symplectic group appears already as 
$\Aut^0 E_n/\Inn E_n$ was probably first done by Isaacs \cite[Section~$4$]{Is}.
Isaacs also computes the character of $\eta$ at odd order elements but we need
to transform the information to a more convenient shape.
}
\end{remar}
The following useful fact demonstrates the multiplicative nature of the Weil
representations. 

\begin{lemma}\label{tp1}
Let $\eta_m$ be a Weil representation of the extrasymplectic group $\ESp(2m,q)$
and let $H$ be an odd order subgroup.

Let $\lambda : \ESp(2m,q) \rightarrow \Sp(2m,q)$
be the natural projection and let $V$ be the natural module for $\Sp(2m,q)$.
Let $h \in \ESp(2m,q)$ be of odd order. Suppose that $\lambda(h)$ preserves
an orthogonal decomposition $V = V_1 \oplus V_2$ and let $m_i = \dim V_i$ for
$i=1, 2$. Then $h = h_1 h_2$ where $h_1, h_2 \in \ESp(2m,q)$, $h_1h_2 = h_2h_1$
and  $\lambda(h_1)$ (respectively, $\lambda(h_2)$) acts trivially on $V_2$ 
(respectively, on $V_1$), and $\eta_m(h) = \eta_{m_1}(h_1) \cdot \eta_{m_2}(h_2)$.
\end{lemma}
\begin{prf}
This is contained in \cite[Lemma 5.5]{Is}. 
\end{prf}

\begin{lemma}\label{tr1} 
Suppose that $n > 1$, let $\eta$ be an irreducible representation of 
$R = \ESp(2n,2)$ as described above, and let $T \subset \Sp(2n,2)$ be a 
maximal cyclic torus of order $2^n+\ep$ where $\ep =1$ or $-1$. Let $T'$ be
any subgroup of $R$ such that $|T'|=|T|$ and $T'E_n/E_n=T$. Then $\chi_{T'} =
\rho_{T'} + \ep \cdot 1_{T'}$ where $\chi$ is the character of $\eta$, that is,
the Weil character of $R$.
\end{lemma}
\begin{prf}
This is a particular case of \cite[Theorem 9.18]{DH}, however, we have to 
refine a few details. Firstly, Theorem 9.18 in \cite{DH} is stated for an 
extraspecial group in place of $E_n$. However, it is known that $T$ is 
contained either in $O^+ (2n,2)$ or in $O^- (2n,2)$ and we can use the
result for extraspecial groups. Secondly, Theorem 9.18 in \cite{DH} claims 
that $\chi_{T'} = \rho_{T'} + \ep \cdot \tau$ where $\tau$ is some linear character
of $T'$. To deduce that in our situation $\tau = 1_{T'}$, observe that $R$
is perfect (unless $n\leq 2$) and hence $\det \eta (t)=1$ for
any $t\in T'$. This is also true for $n=2$ as $\Sp(4,2)$ has a
simple subgroup of index 2, so $T'$ belongs to the derived
subgroup of $R$. As $\det \eta (t) = \tau(t)$, the claim follows.
\end{prf}

\medskip

We fix an embedding $e:\Sp(2m,2^k)\ra
\Sp(2n,2)$ where $n=mk$ and denote by $\ESp(2m,q)$ the preimage of
$\Sp(2m,2^k)$ in $R=R(mk,2)$. Moreover, if $T$ is a maximal torus                                  in $\Sp(2m,2^k)$ then $e(T)$ is a maximal torus in $\Sp(2n,2)$, and
$e(T_1)\times \cdots \times e(T_{k+l})$ is an $e(T)$-decomposition
of $e(T)$. Then Lemmas \ref{tr1} and \ref{tp1} yield the following result.

\begin{propo}\label{pb5}
Let $T$ be a maximal torus in $\Sp(2n,2)$, and let $T=T_1\times \cdots \times                         T_k\times T_{k+1}\times\cdots \times T_{k+l}$ be a $T$-decomposition such that
$|T_i|=2^{n_i}-1$ for $i\leq k$ and $|T_i|=2^{n_i}+1$ for $i>k$.
Let $T',T_{i}'$ be subgroups of $R$ such that $|T'|=|T|$,
$|T_i'|=|T_i|$ for $1\leq i\leq k+l$, and $T'E_n/E_n=T$,
$~T_i'E_n/E_n=T_i$. Let $\chi$ be the character of $\eta$. Then
$$\chi_{T'}=(\rho_{T'_1}+1_{T'_1})\boxtimes\cdots \boxtimes
(\rho_{T'_k}+1_{T'_k})\boxtimes
(\rho_{T'_{k+1}}-1_{T'_{k+1}})\boxtimes\cdots \boxtimes
(\rho_{T'_{k+l}}-1_{T'_{k+l}}).$$ 
Furthermore, this is true for maximal tori in $\Sp(2m,2^k)\subset \Sp(2n,2)$ 
where $n=mk$.
\end{propo}

\begin{remar}
{\rm
It follows that $\chi_{T'}$ is real valued and moreover, that $\chi (g)$ is a
real number for every $g$ of odd order, as the projection of $g$ in $\Sp(2n,2)$
belongs to some maximal torus of $\Sp(2n,2)$.
}
\end{remar}

\subsubsection{The Brauer character of $\sigma_n$}
Recall that $\lambda_1 \ld \lambda_n$ denote the fundamental weights of 
$\mathbf{G}$; for uniformity of some formulas below we set $\lambda_0=0$.
We often use without accurate reference Steinberg's famous theorem saying that
every irreducible representation of $G$ is of shape ${(\phi_\nu)}_{G}$ where 
$\nu$ is a $q$-restricted dominant weight, and conversely ${(\phi_\nu)}_{G}$
is irreducible for every $q$-restricted dominant weight $\nu$ of $\mathbf{G}$.
Recall that a dominant weight $\nu = a_1\lambda_1+\cdots +a_n\lambda_n$ is 
called $q$-restricted if $0\leq a_i\leq q-1$ (here $a_1\ld a_n$ are integers). 
In addition, if $\nu$ is not $2$-restricted then $\phi_\nu$ can be expressed 
as the tensor product of $2$-restricted irreducible representations twisted by
the Frobenius morphism as follows. 
Let $q=2^k$ and let $a_i=\sum _{j=0}^{k-1}2^jb_{ij}$ be the $2$-adic
expansion of $a_i$. Let $\nu_j=\sum_i b_{ij}\lambda_i$. Then
$\phi_\nu =\phi_{\nu_0} \otimes F_0(\phi_{\nu_1}) \otimes \cdots \otimes
F_0^{k-1}(\phi_{\nu_{k-1}})$ where $F_0$ is the standard Frobenius morphism of
$\mathbf{G}$ induced by the mapping $x \mapsto x^2$ for $x \in \mathbf{K}$ 
(so that $F = F_0^k$). In particular, $\phi_{(q-1) \lambda_n}= \phi_{\lambda_n}
\otimes F_0 (\phi_{\lambda_n}) \otimes \cdots \otimes F_0^{k-1}(\phi_{\lambda_n})$;
this fact will be also used without precise reference.

\begin{lemma} \label{nd5}
{\rm \cite[Lemma 1.13]{Z90}} Let $e : \Sp(2n,\mathbf{K}) \ra 
\Sp(2nk,\mathbf{K})$ be the embedding defined by $g \mapsto \diag (g,F_0(g) \ld 
F_0^{k-1}(g))$ for $g \in \Sp(2n,\mathbf{K})$ (this is called a Frobenius 
embedding in {\rm \cite{Z90}}). Then the restriction of $\phi_{\lambda_{nk}}$
to $e(\Sp(2n,\mathbf{K}))$ is irreducible and coincides with $\phi_{\lambda_n}
\otimes F_0(\phi_{\lambda_n})\otimes \cdots \otimes F_0^{k-1}(\phi_{\lambda_n}) 
= \phi_{(q-1)\lambda_n}$. Here, $\phi_{\lambda_{nk}}$ is the irreducible 
representation of $\Sp(2nk,\mathbf{K})$ corresponding to the fundamental weight
$\lambda_{nk}$, while $\phi_{\lambda_n}$ and $\phi_{(q-1)\lambda_n}$ refer to 
the group $\mathbf{G}=\Sp(2n,\mathbf{K})$.
\end{lemma}

\begin{corol}\label{c20}
The restriction
${(\sigma_{nk,2})}_{\Sp(2n,q)}$ is  equivalent to $\sigma_{n,q}$. (He\-re,
$\sigma_{nk,2}$ is the spinor representation of $\Sp(2nk,2)$ and $\Sp(2n,q)$ is
viewed as a subgroup of $\Sp(2nk,2)$ under an embedding obtained by regarding
$\mathbb{F}_q$ as a vector space over $\mathbb{F}_2$.)
\end{corol}
\begin{prf}
Let $\mathbf{V}$ be the natural module for $\Sp(2nk, \mathbf{K})$-module (that
is, the one of highest weight $\lambda_1$).
Then $\mathbf{V}_{\Sp(2n,q)}$ is reducible, in fact $\mathbf{V}_{\Sp(2n,q)} 
\cong V_n \oplus F_0(V_n) \oplus \cdots \oplus F_0^{k-1}(V_n)$ where $V_n$ is
the natural $\Sp(2n,q)$-module. So the result follows from Lemma \ref{nd5}.
\end{prf}

\begin{propo}\label{p88}
{\rm \cite[Theorem 3.10]{Z90}} The Brauer reduction modulo $2$ of $\eta$
is irreducible and equivalent to the inflation of $\sigma_{n,2}$ to $R = R(n,2)$.
\end{propo}

\begin{corol}\label{pc8}
The Brauer character of $\sigma_{n,2}$ is real and coincides with the character of
$\eta$ at elements of odd order. 
\end{corol}

\begin{propo}\label{p99}
Let $q=2^k$. The Brauer reduction modulo $2$ of $\eta$ is irreducible and 
equivalent to the inflation of $\sigma_{n,q}$ to $\ESp(2n,q)$.
\end{propo}
\begin{prf}
This is not explicitly stated in \cite{Z90}, but follows from Corollary 
\ref{c20}. Indeed, by Proposition \ref{p88}, the reduction of $\eta$ modulo $2$
coincides with ${(\sigma_{nk,2})}_{\Sp(2m,q)}$ which is $\sigma_{n,q}$ by 
Corollary \ref{c20}.
\end{prf}

\begin{propo}\label{p55}
Let $T$ be a maximal torus of $\Sp(2n,q)$ and let $T=T_1\times \cdots \times 
T_k\times T_{k+1}\times\cdots \times T_{k+l}$ be a $T$-decomposition such that 
$|T_i|=q^{n_i}-1$ for $i\leq k$ and $|T_i|=q^{n_i}+1$ for $i>k$. Let $\omega$ 
be the Brauer character of $\sigma_{n,q}$. Then 
$$\omega_T = (\rho_{T_1}+1_{T_1}) \boxtimes \cdots \boxtimes 
(\rho_{T_k}+1_{T_k}) \boxtimes (\rho_{T_{k+1}}-1_{T_{k+1}}) \boxtimes \cdots
\boxtimes (\rho_{T_{k+l}}-1_{T_{k+l}}).$$
\end{propo}
\begin{prf}
This follows from Propositions \ref{pb5} and \ref{p99}.
\end{prf}

\begin{corol} \label{c75}
Let $g \in \Sp(2n,q)$ be an odd order element. Then $\omega(g)^2
= q^{N(V;g)}$ where $V$ is the natural $\Sp(2n,q)$-module and $N(V;g)$ the 
dimension of the $1$-eigenspace of~$g$ on~$V$.
\end{corol}
\begin{prf}
This can be deduced from Lemma \ref{tr1} but is also available in Isaacs 
\cite[Theorem 3.5]{Is}.
\end{prf}

\subsection{Multiplicities in $\om_T$}

We return to the general situation. Namely, $G = \Sp(2n,q)$ or $U(d,q)$
with~$q$ arbitrary. If $G = \Sp(2n,q)$ and $q$ is even, we let $\omega$
denote the (Brauer) character of~$G$ of the representation~$\sigma_{n,q}$
as in Subsection~\ref{ModularCase}. Otherwise,~$\omega$ denotes the character
of a Weil representation of~$G$ as introduced in Subsection~\ref{OrdinaryCase}.
We let~$T$ be a maximal torus of~$G$ and consider a $T$-decomposition
$$T = T_1 \times \cdots \times T_k \times T_{k+1} \times \cdots \times T_{k+l}$$
as in Subsection~\ref{TDecomposition}.

Let $\theta \in \Irr(T)$. 
Then $\theta = \theta_1 \boxtimes \cdots \boxtimes \theta_{k+l}$
with unique $\theta_i \in \Irr(T_i)$, $1 \leq i \leq k + l$. If $k = 1$ and
$l=0$, then we see from Lemma~\ref{g1}(a) and Proposition~\ref{p55}, 
respectively, that the multiplicity of every $\theta \in \Irr(T)$ in 
$\omega_T$ equals~$1$, except for the character
$1^-_T$, which has multiplicity~$2$. (Recall our convention that $1_T^-$ stands
for $1_T$ if~$T$ has odd order, i.e., if~$q$ is even.)
Similarly, if $k=0$ and $l=1$, then
Lemma~\ref{g1}(b) respectively Proposition~\ref{p55} implies that the 
multiplicity of every
$\theta \in \Irr(T)$ in $\omega_T$ equals~$1$, except for the character
$1^-_T$, which has multiplicity~$0$. In general, let $k(\theta)$ be the
number of $i \leq k$ such that  $\theta_i=1^-_{T_i}$. It follows from Lemma
\ref{g1}(c) that the multiplicity of $\theta \in \Irr(T)$ in $\omega_T$ equals
$2^{k(\theta)}$, unless there is $j$ such that $\theta_{k+j} =1^-_{T_{k+j}}$,
in which case $\theta$ does not occur in $\omega_T$. Thus we have proved the
following.

\begin{lemma}\label{g2}
Let $\theta = \theta_1 \boxtimes \cdots \boxtimes \theta_{k+l}$ be an irreducible
character of $T = T_1 \times \cdots \times  T_{k+l}$.

$(1)$ If $\theta_{k+j} =1^-_{T_{k+j}}$ for some $j>0$ then $\theta$ does not
occur as an irreducible constituent of $\omega_T$ (that is, $(\om_T,\theta )=0$).

$(2)$ Suppose that $\theta_{k+j} \neq 1^-_{T_{k+j}}$ for every $j=1, \ldots , l$.
Let $k(\theta)$ be the number of $0 \leq i \leq k$ such that $\theta_i=1^-_{T_i}$.
Then $(\om_T,\theta )= 2^{k(\theta)}$.

$(3)$ Suppose that $\theta_{i} \neq 1^-_{T_{i}}$ for every $1 \leq i \leq
k + l$. Then $(\om_T,\theta )= 1$.
\end{lemma}
\noindent
Note that the statements above remain true in case~$G$ is a group
of characteristic~$2$, if~$1_{T_i}^-$  is replaced by $1_{T_i}$
throughout, in consistency with our convention.

\section{The product $\om \cdot {\rm \St}$}

In this section we prove Theorem~\ref{Thm1Intro} for the symplectic and 
unitary groups.

Let $(\mathbf{G}, F)$ be a unitary group as in~\ref{UnitaryGroups} or a 
symplectic group as in~\ref{SymplecticGroups}. 
We denote by $(\mathbf{G}^*, F^*)$ a reductive group dual to $(\mathbf{G}, F)$.
Let $S^*$ denote the set of $G^*$-conjugacy classes of semisimple elements
of~$G^*$. We write $(s^*)$ for the element of ${S}^*$ containing $s^* \in G^*$.
For each semisimple $s^* \in G^*$ we choose a set ${\kappa}(s^*)$ of
representatives for the $G$-orbits in 
$$\{ (\mathbf{T},\theta) \in
\mathcal{T}(\mathbf{G}) \mid (\mathbf{T},\theta) \text{\rm\ is dual to }
(\mathbf{T}^*,s^*) \in \mathcal{S}(\mathbf{G}^*) \}$$
(see~(\ref{TGandSG})).

By $\St = \St_{G}$ we denote the Steinberg character of $G = \mathbf{G}^F$ and
by~$\omega$ the class function introduced in Section~\ref{WeilCharacters}. Then
$\omega \cdot \St$ is a character of~$G$ vanishing on all $p$-singular elements.
 It is known that every such class function is uniform, that is, a linear 
combination of characters $R_{\mathbf{T},\theta}$ (see \cite[page 89]{DM}).

The argument in \cite[p.~242]{C} shows that
\begin{equation}
\label{EqOne}
\omega \cdot \St =
\sum_{(s^*) \in {S}^*}~~~\sum _{(\mathbf{T},\theta)\in {\kappa}(s^*) }
\frac{\ep_\mathbf{G}\ep_\mathbf{T}
(\omega_T,\theta)}{|W(\mathbf{T})^F_\theta|} R_{\mathbf{T},\theta}.
\end{equation}
For each $(s^*) \in {S}^*$ consider the partial sum
\begin{equation}
\label{EqThree}
\pi_{s^*} := \sum_{(\mathbf{T},\theta)\in \kappa(s^*) }
\frac{\ep_\mathbf{G}\ep_\mathbf{T}(\om_T,\theta)}{|W(\mathbf{T})^F_\theta|}
R_{\mathbf{T},\theta},
\end{equation}
as well as the class function
\begin{equation}
\label{EqFour}
\rho_{s^*}= \sum_{(\mathbf{T},\theta)\in {\kappa}(s^*)}
\frac{\ep_\mathbf{G}\ep_\mathbf{T}}{|W(\mathbf{T})^F_{\theta}|}
R_{\mathbf{T},\theta}.
\end{equation}

\begin{lemma}\label{RhosDiMi}
Let $\chi_{(s^*)}$ be the class function introduced by Digne and Michel
in~\cite[Definition 14.10]{DM}. Then $\chi_{(s^*)} = \rho_{s^*}$ if
$C_{\mathbf{G}^*}(s^*)$ is connected.
\end{lemma}
\begin{prf}
In the notation of~\cite{DM}, 
$$\chi_{(s^*)} := |W^\circ(s^*)|^{-1}\sum_{w \in W^\circ(s^*)} 
\ep_\mathbf{G}\ep_{\mathbf{T}_w^*} R_{\mathbf{T}_w^*}(s^*).$$
We begin by explaining this notation. Firstly, $W^\circ(s^*)$ is the Weyl 
group of $C^\circ_{\mathbf{G}^*}(s^*)$, the connected component of 
$C_{\mathbf{G}^*}(s^*)$. Since $C_{\mathbf{G}^*}(s^*)$ is connected, we have
$C^\circ_{\mathbf{G}^*}(s^*) = C_{\mathbf{G}^*}(s^*)$ and hence $W^\circ(s^*) 
= W(s^*)$ (see \cite[Remark~$2.4$]{DM}). Secondly, $R_{\mathbf{T}_w^*}(s^*)$ 
denotes a Deligne-Lusztig character of~$G$ of the form 
$R_{\mathbf{T},\vartheta}$, where $(\mathbf{T},\vartheta) \in
\mathcal{T}(\mathbf{G})$ is dual to $(\mathbf{T}_w^*,s^*) \in 
\mathcal{S}(\mathbf{G}^*)$, and where $\mathbf{T}_w^*$ is obtained from the
reference torus of $C_{\mathbf{G}^*}(s^*)$ by twisting with~$w$ (cf.\ 
Subsection~\ref{ToriClassification}).

Let $\kappa^*(s^*)$ denote a set of representatives for the $F^*$-conjugacy 
classes of $W(s^*)$. Then, again by the results summarized in 
Subsection~\ref{ToriClassification}, we have
$$\chi_{(s^*)} = \sum_{w \in \kappa^*(s^*)}
\frac{\ep_\mathbf{G}\ep_{\mathbf{T}^*_{w}}}{|C_{W(s^*),F^*}(w)|} 
R_{\mathbf{T}^*_{w}}(s^*).$$
Every element of $\kappa(s^*)$ is dual (in the sense of~(\ref{TGandSG})) 
to a pair $({\mathbf{T}'}^*, s^*) \in \mathcal{S}(\mathbf{G}^*)$; since 
$s^* \in {\mathbf{T}'}^*$, we have in fact $({\mathbf{T}'}^*, s^*) \in 
\mathcal{S}(C_{\mathbf{G}^*}(s^*))$. Two such pairs are conjugate in $G^*$ if
and only if they are conjugate in $C_{\mathbf{G}^*}(s^*)^{F^*}$. Thus there is
a bijection $\kappa^*(s^*) \rightarrow \kappa(s^*)$ such that
$(\mathbf{T},\vartheta) \in \kappa(s^*)$ is dual to $(\mathbf{T}^*_{w},s^*)$
if $w \in \kappa^*(s^*)$ is mapped to $(\mathbf{T},\vartheta)$.
By Lemma~\ref{SameOrders} we have $|W(\mathbf{T})^F_\theta| = 
|W(\mathbf{T}^*_{w})^{F^*}_{s^*}|$ for pairs corresponding this way.

Note that $W(\mathbf{T}^*_{w})^{F^*}_{s^*} 
= N_{C_{\mathbf{G}(s^*)}}({\mathbf{T}^*_{w}})^{F^*}/{\mathbf{T}^*_{w}}^{F^*} 
\cong C_{W(s^*),F^*}(w)$, the latter by \cite[Proposition~$3.3.6$]{C}, applied to 
$C_{\mathbf{G}^*}(s^*)$. This completes the proof.
\end{prf}

The above result does not hold if $C_{\mathbf{G}^*}(s^*)$ is not connected. 
Consider, for example, the case $\mathbf{G} = \Sp(2,\mathbf{K}) \cong \SL(2, 
\mathbf{K})$, where~$q$ is odd. There is an involution $s^* \in G^* = \SO(3,q)
 \cong \PGL(2,q)$ whose centralizer is equal to 
$N_{\mathbf{G}^*}(\mathbf{T}_0^*)$. If $(\mathbf{T_0}, \theta)$ is dual to 
$(\mathbf{T}_0^*,s^*)$, then $\theta = 1_{T_0}^-$, and 
$|W(\mathbf{T}_0)^F_{\theta}| = 2$. Since $|W^\circ(s^*)| = 1$, we have 
$\chi_{(s^*)} = 2\rho_{s^*}$.

For the sake of a uniform notation, we introduce a basis $v_1^*$,
$v_2^*, \ldots, v_n^*$, [$v_0^*$,] ${v_n^*}', \ldots, {v_2^*}'$,
${v_1^*}'$ of the vector space $V^*$ (where $v_0^*$ is not present
if $d = \dim V^*$ is even), such that $v_1^*, v_2^*, \ldots,
v_n^*$ and ${v_n^*}', \ldots, {v_2^*}', {v_1^*}'$ span maximal
isotropic subspaces of $V^*$ and the hermitean or orthogonal form
takes value~$1$ on the pairs $v_i^*, {v_i^*}'$, $1 \leq i \leq n$,
and $v_0^*$, if present, has norm~$1$. (Thus in the orthogonal
case we have just ``starred'' the basis
from~\ref{SymplecticGroups}.)

\begin{lemma}\label{kn1}
Let $s^*$ be a semisimple element of $G^*$ without eigenvalue $(-1)^q$ on
$V^*$. Then $\pi_{s^*} = \rho_{s^*} \in \Irr(G)$.
\end{lemma}
\begin{prf}
We have $(\om_T,\theta) = 1$ for all $(\mathbf{T}, \theta) \in {\kappa}(s^*)$
by Lemmas \ref{g2}(3) and~\ref{SameOrders}(b). Hence the expression for
$\pi_{s^*}$ coincides with that for $\rho_{s^*}$.

Now $C_{\mathbf{G}^*}(s^*)$ is connected since~$s^*$ does not have
eigenvalue~$(-1)^q$. (If $\mathbf{G}^* = \GL_n(\mathbf{K})$ the
centralizer of every semisimple element is connected. In the other
case, the result can be derived from \cite[Theorem 3.5.3]{C}.) The
irreducibility of $\rho_{s^*}$ follows from Lusztig's results in
\cite{Lu} (see also \cite[$14.40$, $14.43$, $14.48$]{DM} in connection with
Lemma~\ref{RhosDiMi}).
\end{prf}

Our goal now is to determine the class functions $\pi_{s^*}$ in case~$s^*$ has
eigenvalue $(-1)^q$ on~$V^*$.

\begin{lemma}
\label{ConjugateToStandard}
Let $s^* \in G^*$ be a semisimple element which has eigenvalue $(-1)^q$ on~$V^*$
and suppose that $\pi_{s^*} \neq 0$.
Then~$s^*$ is conjugate in $\mathbf{G}^*$ to an element whose $(-1)^q$-eigenspace
on~$V^*$ equals $\langle v_1^*, \ldots , v_m^*, {v_m^*}', \ldots , {v_1^*}' \rangle$
for some $1 \leq m \leq n$.
\end{lemma}
\begin{prf}
Denote by $V^*_-$ the $(-1)^q$-eigenspace of~$s^*$, and by ${(V^*_-)}^\perp$
its orthogonal complement. Every element of~$G^*$ commuting with~$s^*$ fixes
$V^*_-$ and ${(V^*_-)}^\perp$, and so every maximal torus $T^*$ of $G^*$
containing~$s^*$ yields a $T^*$-decomposition of~$V^*$ compatible
with the direct sum $V^* = V^*_- \oplus {(V^*_-)}^\perp$.

Let $T^*$ be a maximal torus of $G^*$ with $s^* \in T^*$ and let
$V^* = V^*_0 \oplus V^*_1 \oplus \cdots \oplus V^*_k \oplus V^*_{k+1}
\oplus \cdots \oplus V^*_{k+l}$ be such a compatible $T^*$-decomposition.

Suppose first that $(\mathbf{G}^*,F^*)$ is unitary, and that $\dim
V^*_-$ is odd. Then $V^*_j \subseteq V^*_-$ for some $j > k$ by
Lemma~\ref{n1}(1). Now let $(\mathbf{G}^*,F^*)$ be orthogonal.
Then $V^*_-$ has even dimension $2m$. Suppose that the Witt index
of $V^*_-$ is smaller than~$m$. Then, again, $V^*_j \subseteq
V^*_-$ for some $j > k$ by Lemma~\ref{n1}(2). It follows that in
the decomposition of the corresponding pair $(T,\theta)$, we have
$\theta_j = 1^-_{T_j}$ by Lemma~\ref{SameOrders}. By
Lemma~\ref{g2}(1) this implies that $(\omega_T,\theta) = 0$. Thus
$\pi_{s^*} = 0$ contrary to our assumption.

Hence $\dim V^*_- = 2m$ is even, and the Witt index of $V^*_-$ equals~$m$
in the orthogonal case. By Witt's theorem we may assume that
$V^*_- = \langle v_1^*, \ldots , v_m^*, {v_m^*}', \ldots , {v_1^*}' \rangle$.
\end{prf}

For $1 \leq m \leq n$ write $${\mathbf{V}^{(m)}}^* := \langle
v_1^*, \ldots , v_m^*, {v_m^*}', \ldots , {v_1^*}'
\rangle_{\mathbf{K}},$$ and $${\mathbf{V}^{(m')}}^* := \langle
v_{m+1}^*, \ldots , v_n^*, [v_0^*,] {v_n^*}', \ldots ,
{v_{m+1}^*}' \rangle_{\mathbf{K}},$$ where the notation $[v_0^*]$
indicates that $v_0^*$ is to be omitted if $\dim \mathbf{V}^*$ is
even. As usual we denote the sets of rational points of these
vector spaces by ${V^{(m)}}^*$ and ${V^{(m')}}^*$, respectively.
Then ${V^{(m')}}^*$ is the orthogonal complement of ${V^{(m)}}^*$.

Let ${\mathbf{G}^{(m)}}^*$ denote the subgroup of $\mathbf{G}^*$ fixing
${\mathbf{V}^{(m)}}^*$ and acting as the identity on ${\mathbf{V}^{(m')}}^*$,
and let ${\mathbf{G}^{(m')}}^*$ be defined similarly. Then ${\mathbf{G}^{(m)}}^*
\times {\mathbf{G}^{(m')}}^* \leq \mathbf{G}^*$ is the identity component of
the stabilizer in $\mathbf{G}^*$ of the direct sum decomposition
$\mathbf{V}^* = {\mathbf{V}^{(m)}}^* \oplus {\mathbf{V}^{(m')}}^*$.

\begin{lemma}
\label{TorusInLevi}
Fix $1 \leq m \leq n$ and let $s^* \in G^*$ be a semisimple element whose
$(-1)^q$-eigenspace on~$V^*$ equals ${V^{(m)}}^* = \langle v_1^*, \ldots ,
v_m^*, {v_m^*}', \ldots , {v_1^*}' \rangle$.

Let $\mathbf{T}^*$ be an $F^*$-stable maximal torus of $\mathbf{G}^*$
containing~$s^*$ and let $(\mathbf{T},\theta) \in \mathcal{T}(\mathbf{G})$
be dual to $(\mathbf{T}^*,s^*)$. If $(\omega_T,\theta) \neq 0$, then
$\mathbf{T}^*$ is conjugate
in $G^*$ to a torus fixing $\langle v_1^*, \ldots , v_m^* \rangle_{\mathbf{K}}$.
\end{lemma}
\begin{prf}
Clearly, $\mathbf{T}^*$ fixes the $(-1)^q$-eigenspace ${\mathbf{V}^{(m)}}^*$ 
of $s^*$ and its orthogonal complement ${\mathbf{V}^{(m')}}^*$.
Consider a $T^*$-decomposition $V^* = V_0^* \oplus V_1^* \oplus \cdots \oplus V_k^*
\oplus V_{k+1}^* \oplus \cdots \oplus V_{k+l}^*$ of $V^*$ compatible with the orthogonal
decomposition $V^* = {V^{(m)}}^* \oplus {V^{(m')}}^*$. If $V_j^* \leq {V^{(m)}}^*$
for some $j > k$, then, in the decomposition of the corresponding pair $(T,\theta)$, we
have $\theta_j = 1^-_{T_j}$ by Lemma~\ref{SameOrders}, and so $(\omega_T,\theta) = 0$
by Lemma~\ref{g2}(1). Thus our assumption implies that ${V^{(m)}}^*$ is a direct
sum of some $V_j^*$s with $1 \leq j \leq k$, and so $T^*$ fixes a maximal
singular subspace of ${V^{(m)}}^*$.
By conjugating $\mathbf{T}^*$ by an element of $G^*$, we may assume that
$T^*$ fixes $\langle v_1^*, \ldots , v_m^* \rangle$.

If $\dim \mathbf{V}^*$ is odd, we may also assume that $\mathbf{T}^*$ fixes $v_0^*$,
by conjugating $\mathbf{T}^*$ with a suitable element of ${G^{(m')}}^*$. It follows
that $\mathbf{T}^*$ fixes the space $\langle v_0^* \rangle_{\mathbf{K}} \oplus
{\mathbf{V}^{(m)}}^*$ in this case. We may thus assume that $m = n$. Using the
classification of the maximal tori in Sections~\ref{UnitaryGroups}
and~\ref{SymplecticGroups},
we see that $\mathbf{T}^*$ is conjugate in $G^*$ to a maximal $F^*$-stable torus
fixing $\langle v_1^*, \ldots , v_n^* \rangle_{\mathbf{K}}$. (If this were not the
case, then $T^*$ would have an irreducible direct summand different from
$\langle v_0^* \rangle$ in a $T^*$-decomposition of~$V^*$.)
\end{prf}

Thus we may assume that every pair $(\mathbf{T},\theta)$ which contributes
a non-zero summand to the sum~(\ref{EqThree}) is dual to a pair
$(\mathbf{T}^*,s^*)$ such that $\mathbf{T}^*$ fixes $\langle v_1^*,
\ldots , v_m^* \rangle_{\mathbf{K}}$ for some $1 \leq m \leq n$. In other
words, $\mathbf{T}^*$ lies in the standard (split) Levi subgroup
${\mathbf{L}^{(m)}}^* \times {\mathbf{G}^{(m')}}^*$ of $\mathbf{G}^*$ fixing
$\langle v_1^*, \ldots , v_m^* \rangle_{\mathbf{K}}$. Here, ${\mathbf{L}^{(m)}}^*$
denotes the standard Levi subgroup of ${\mathbf{G}^{(m)}}^*$ fixing
$\langle v_1^*, \ldots , v_m^* \rangle_{\mathbf{K}}$. Moreover, two such
tori are conjugate in $G^*$ if and only if they are conjugate in ${L^{(m)}}^*
\times {G^{(m')}}^*$.

We now fix $1 \leq m \leq n$, an element $s^* \in G^*$ whose
$(-1)^q$-eigenspace on $V^*$ equals ${V^{(m)}}^*$, and a maximal torus
$\mathbf{T}^* \leq {\mathbf{L}^{(m)}}^* \times {\mathbf{G}^{(m')}}^*$
containing~$s^*$. Let $(\mathbf{T},\theta) \in \mathcal{T}(\mathbf{G})$
be a pair dual to $(\mathbf{T}^*,s^*)$.
Since duality behaves well with respect to split Levi subgroups,
we may assume that $\mathbf{T} \leq \mathbf{L}^{(m)} \times \mathbf{G}^{(m')}$,
the standard Levi subgroup of $\mathbf{G}$ fixing the isotropic subspace
$\langle v_1, \ldots , v_m \rangle_{\mathbf{K}}$ of $\mathbf{V}$.

We have $\mathbf{L}^{(m)} \cong \GL(m,\mathbf{K})$ (acting on $\langle v_1, \ldots ,
v_m \rangle_{\mathbf{K}}$). Furthermore, we may assume that $\mathbf{T} =
\mathbf{T}^{(m)} \times
\mathbf{T}^{(m')}$ with $F$-stable maximal tori of $\mathbf{L}^{(m)}$ and of
$\mathbf{G}^{(m')}$, respectively, and we have a corresponding decomposition
$\theta = \theta^{(m)} \boxtimes \theta^{(m')}$.

To simplify notation, we put $\mathbf{L} := \mathbf{L}^{(m)}$, $\mathbf{G}' :=
\mathbf{G}^{(m')}$, $\mathbf{S} := \mathbf{T}^{(m)}$, $\mathbf{T}' :=
\mathbf{T}^{(m')}$, $\sigma := \theta^{(m)}$ and $\theta' := \theta^{(m')}$.
Then $\mathbf{T} = \mathbf{S} \times \mathbf{T}'$ and $\theta = \sigma \boxtimes
\theta'$.

\begin{lemma}\label{kn2}
With the above notation we have:
$$ \frac{ (\om_T,\theta ) }{ |W(\mathbf{T})^F_\theta|}
=\frac{1}{|W_{\mathbf{L}}(\mathbf{S})^F_{\sigma}|} \cdot
\frac{ 1 }{|W_{\mathbf{G}'}(\mathbf{T}')^F_{\theta'}|}.$$
\end{lemma}
\begin{prf}
Let $\mathbf{H}$ denote the subgroup of $\mathbf{G}$ fixing $\mathbf{V}^{(m)}$
and acting as the identity on its complement $\mathbf{V}^{(m')}$. Then $\mathbf{H}$
is a general linear or symplectic group of dimension $2m$ over~$\mathbf{K}$.
Moreover,~$S$ is a neutral maximal torus of~$H$.

Using the multiplicity of the Weil representation (see
Section~\ref{WeilCharacters}) and Lemma \ref{WeylGroupFactoriazation}
we find
$$\frac{(\om_T,\theta)}{|W(\mathbf{T})^F_\theta|} =
\frac{(\om^{(H)}_{S},\sigma)}{|W_{\mathbf{H}}(\mathbf{S})^F_{\sigma}|}
\cdot \frac{(\om^{(G')}_{T'},\theta')}
{|W_{\mathbf{G}'}(\mathbf{T}')^F_{\theta'}|}.$$
Now $(\om^{(G')}_{T'},\theta') = 1 $ by Lemma~\ref{g2}(3).
The claim follows as long as we can show that
$$
\frac{(\om^{(H)}_{S},\sigma)}{|W_{\mathbf{H}}(\mathbf{S})^F_{\sigma}|}
= \frac{1}{|W_{\mathbf{L}}(\mathbf{S})^F_{\sigma}|}.$$
Let $(1^{l_1}, 2^{l_2}, \dots , m^{l_m})$ be the partition of~$m$ defining
the neutral maximal torus $\mathbf{S}$ of $\mathbf{H}$ (see Lemma
\ref{n1}(3)). By Lemma \ref{g2}(2), $(\om^{(H)}_{S},\sigma) =
2^{k(\sigma)}=2^{l_1+\cdots + l_m}$. By Lemma \ref{n2},
$|W_{\mathbf{H}}(\mathbf{S})^F|=
2^{l_1+\cdots +l_m}|W_{\mathbf{L}}(\mathbf{S})^F|$, proving the
desired result.
\end{prf}

We return to the computation of $\pi_{s^*}$, with $s^*$ as above. Write
$s^* = (-1)^q \times {s^*}'$ with $(-1)^q \in {L^{(m)}}^*$ and ${s^*}' \in
{G^{(m')}}^*$. By the considerations above, if $(\mathbf{T},\theta) \in
\mathcal{T}(\mathbf{G})$ is dual to $(\mathbf{T}^*, s^*)$, and if
$(\om_T,\theta) \neq 0$, then we may assume that there is a factorisation
$$(\mathbf{T},\theta) = (\mathbf{S} \times \mathbf{T}', \sigma
\boxtimes \theta'),$$
in such a way that $(\mathbf{S},\sigma) \in \mathcal{T}(\mathbf{S})$
is dual to $(\mathbf{S}^*,(-1)^q) \in \mathcal{S}(\mathbf{L}^*)$ and
$(\mathbf{T}',\theta') \in \mathcal{T}(\mathbf{G}')$
is dual to $({\mathbf{T}'}^*,{s^*}') \in \mathcal{S}({\mathbf{G}'}^*)$.
Thus we may restrict summation in~(\ref{EqThree}) to
$\kappa_{\mathbf{L} \times \mathbf{G}'}(s^*) = {\kappa}_{\mathbf{L}}((-1)^q)
\times {\kappa}_{\mathbf{G}'}({s^*}')$, with the obvious interpretation
of ${\kappa}_{\mathbf{L}}$ and ${\kappa}_{\mathbf{G}'}$.

Let~$\mathbf{P}$ denote the the standard parabolic subgroup of $\mathbf{G}$
fixing the isotropic subspace $\langle v_1, \ldots , v_m \rangle_{\mathbf{K}}$
of $\mathbf{V}$.
By \cite[Proposition 7.4.4]{C}, $R^\mathbf{G}_{\mathbf{T},\theta} =
\left(\infl_P\left({R^{\mathbf{L} \times \mathbf{G}'}_{\mathbf{T},\theta}}\right)\right)^G$
where $\infl_P\left( \psi \right)$ denotes the inflation of the class
function~$\psi$ of $L \times G'$ to~$P$ via the homomorphism
$P \rightarrow L \times G'$.

\begin{lemma}\label{nf7}
Let the notation be as above. Since $s^* \leq \mathbf{L}^* \times {\mathbf{G}'}^*$,
we have a class function $\rho_{s^*}^{(L \times G')}$ of $L \times G'$ defined
analogously to $\rho_{s^*}$ for~$G$.
With this notation we have $\pi _{s^*}=
\left(\infl_P\left(\rho_{s^*}^{(L \times G')}\right)\right)^G$.
In addition, $\rho_{s^*}^{(L \times G')} = \St_{L}^- \boxtimes \rho^{(G')}_{{s^*}'}$,
where $\St_{L}^- = 1^-_{L} \cdot \St_{L}$.
\end{lemma}
\begin{prf}
We have
\begin{eqnarray*}
\pi_{s^*} & = & \sum_{(\mathbf{T},\theta) \in {\kappa}(s^*) }
\frac{\ep_\mathbf{G}\ep_\mathbf{T}
(\om_T,\theta)}{|W(\mathbf{T})^F_\theta|} R^\mathbf{G}_{\mathbf{T},\theta} \\
 & =  & \left(\infl_P\left(
\sum_{(\mathbf{T},\theta) \in \kappa_{\mathbf{L} \times \mathbf{G}'}(s^*) }
\frac{\ep_\mathbf{G}\ep_\mathbf{T}
(\om_T,\theta)}{|W(\mathbf{T})^F_\theta|}
R^{\mathbf{L} \times \mathbf{G}'}_{\mathbf{T},\theta}\right)\right)^G.
\end{eqnarray*}
By Lemma \ref{kn2} and the discussion above, we find
$$
\sum_{(\mathbf{T},\theta) \in \kappa_{\mathbf{L} \times \mathbf{G}'}(s^*) }
\frac{\ep_\mathbf{G}\ep_\mathbf{T}
(\om_T,\theta)}{|W(\mathbf{T})^F_\theta|}
R^{\mathbf{L} \times \mathbf{G}'}_{\mathbf{T},\theta} =
$$
$$
\sum_{(\mathbf{S},\sigma) \in {\kappa}_{\mathbf{L}}((-1)^q)}
\sum_{(\mathbf{T}',\theta') \in {\kappa}_{\mathbf{G}'}({s^*}')}
\frac{\ep_{\mathbf{L}} \ep_{\mathbf{S}}}
{|W_{\mathbf{L}}(\mathbf{S})^F_{\sigma}| }\cdot
\frac{\ep_{\mathbf{G}'} \ep_{\mathbf{T}'}}
{|W_{\mathbf{G}'}(\mathbf{T}')^F_{\theta'}|}
R^{\mathbf{L}\times\mathbf{G}'}_{\mathbf{T},\sigma \boxtimes \theta'}. $$
Observe that $R^{\mathbf{L}\times\mathbf{G}'}_{\mathbf{T},\sigma \boxtimes \theta'} =
R_{\mathbf{S},\sigma}^{\mathbf{L}} \boxtimes
R_{\mathbf{T}',\theta'}^{\mathbf{G}'}$.
Therefore, the right hand side of the above expression equals the product
$$
\left( \sum_{(\mathbf{S},\sigma) \in {\kappa}_{\mathbf{L}}((-1)^q)}
\frac{\ep_{\mathbf{L}} \ep_{\mathbf{S}}}
{|W_{\mathbf{L}}(\mathbf{S})^F_{\sigma}| }
R_{\mathbf{S},\sigma}^{\mathbf{L}}\right)
\boxtimes
\left(
\sum_{(\mathbf{T}',\theta') \in {\kappa}_{\mathbf{G}'}({s^*}')}
\frac{\ep_{\mathbf{G}'} \ep_{\mathbf{T}'}}
{|W_{\mathbf{G}'}(\mathbf{T}')^F_{\theta'}|}
R_{\mathbf{T}',\theta'}^{\mathbf{G}'}\right).
$$

We have $\sigma = 1_S^-$ for all pairs $(\mathbf{S},\sigma)$ occurring in
the above sum. Hence $R_{\mathbf{S},\sigma}^{\mathbf{L}} = 1_{L}^- \cdot
R_{\mathbf{S},1_{S}}^{\mathbf{L}}$ and $W_{\mathbf{L}}(\mathbf{S})^F_{\sigma} =
W_{\mathbf{L}}(\mathbf{S})^F_{1_{S}}$ for all such pairs. It follows that
the first of these factors equals $\St_{L}^-$ (see \cite[Corollary 7.6.6]{C}),
while the second one, by definition, is equal to $\rho^{(G')}_{{s'}^*}$.
Note that the latter is an irreducible character by~Lemma~\ref{kn1}.
\end{prf}

\smallskip

\textbf{Proof of Theorem~\ref{Thm1Intro} (Part~I)}.
Set $\gamma = \sum _{(s^*) \in {S}^*}\rho _{s^*}$. If $Z(\mathbf{G})$ is 
connected then 
$\gamma$ is  known to coincide with the Gelfand-Graev character of~$G$. Denote
by $\gamma'$ the ``truncated'' character obtained from $\gamma$ by removing all
$\rho_{s^*}$ with $s^*$ having eigenvalue $(-1)^q$. 

Now~(\ref{EqOne}),~(\ref{EqFour}), Lemma~\ref{kn1}, and Lemma~\ref{nf7} yield 
a proof of Theorem~\ref{Thm1Intro} for the symplectic and unitary groups.

\medskip

\textbf{Proof of Corollary~\ref{MultiplicityFree}}.
If $(s_1^*)$ and $(s_2^*)$ are distinct elements of~$S^*$, the constituents
of $\pi_{s_2^*}$ and $\pi_{s_2^*}$ lie in distinct Lusztig series of characters.
Hence it suffices to show that $\pi_{s^*}$ is multiplicity free, if $s^*$ has 
a $2m$-dimensional $(-1)^q$-eigenspace for some $1 \leq m \leq n$.

Lemma~\ref{nf7} shows that $\pi_{s^*} = \left(\infl_P(\St_{L}^- \boxtimes
\rho^{(G')}_{{s^*}'})\right)^G$.
We may use Harish-Chandra theory to see that this Harish-Chandra induced
character is multiplicity free. If~$D$ denotes the maximally split torus of~$L$,
then clearly $\St_{L}^-$ lies in the $(D,1_D^-)$ Harish-Chandra series of~$L$.
Let~$M'$ be a Levi subgroup of $G'$ and $\tau'$ an irreducible cuspidal character
of~$M'$ such that $\rho^{(G')}_{{s^*}'}$ lies in the $(M',\tau')$ Harish-Chandra
series of~$G'$.
Then all constituents of $\pi_{s^*}$ and the irreducible character
$\St_{L}^- \boxtimes \rho^{(G')}_{{s^*}'}$ lie in the
$(D \times M', 1_{D}^- \boxtimes \tau')$ Harish-Chandra series of~$G$ and of
$L \times G'$, respectively.

Now $W_G(D \times M', 1_{D}^- \boxtimes \tau') = W_H(D,1_D^-) \times
W_{G'}(M',\tau')$, by Lemma~\ref{WeylGroupFactoriazation}.
Here,~$H$ has the same meaning as in the proof of Lemma~\ref{kn2}.
Clearly, $W_{L\times G'}(D \times M', 1_{D}^- \boxtimes \tau')
= W_L(D,1_D^-) \times W_{G'}(M',\tau')$. Now $W_H$ and $W_L$ are the
Weyl groups of $H$ and $L$, respectively, the former of type~$B_m$,
the latter its parabolic subgroup, of type~$A_{m-1}$, obtained by
deleting the outer node on the double bond of the Dynkin diagram
for~$W_H$. Via Harish-Chandra theory, the character
$\St_{L}^- \boxtimes \rho^{(G')}_{{s^*}'}$ corresponds to a character
$\mbox{\rm sgn} \boxtimes \lambda'$, where $\mbox{\rm sgn}$ is the
sign character of the symmetric group $W_L \cong S_m$ and $\lambda'$
is some irreducible character of $W' := W_{G'}(M',\tau')$.

By a result of Curtis (see \cite[Theorem (70.24)]{CRII}), the multiplicities
of the irreducible constituents of $\left(\infl_P(\St_{L}^- \boxtimes
\rho^{(G')}_{{s^*}'})\right)^G$ can be computed from the multiplicities
of the induced character $(\mbox{\rm sgn})_{W_L}^{W_H}$. The latter is
the sum of all irreducible characters of $W_H$ which are labelled by
bipartitions of $m$ whose parts are all equal to $1$. (This fact can
be derived from a special case of the Littlewood-Richardson rule; see,
e.g., \cite[Lemma 6.1.4]{GP}.) This completes the proof.

\section{The Weil representation of the general linear group}
\label{WeilInGeneralLinear}

Here, we consider the tensor product of the Weil representation of
the general linear group with its Steinberg representation.

G{\'e}rardin defined the Weil representation of $G := \GL(n,q)$ as the 
permutation representation of $G$ on the vectors of the underlying vector
space (see \cite[Corollary 1.4]{G}). According to 
Definition~\ref{DefinitionOmega}, let us write $\hat{\omega}$ for this 
permutation character and put $\omega := 1_G^-\cdot \hat{\omega}$. 
We will compute $\hat{\omega} \cdot \St$, from which the desired result
follows.

In order to proceed, we describe the stabilizer in $G$ of a non-zero
vector, and its characters. For inductive reasons, we treat $n$, the
dimension of the underlying vector space, as a parameter. In particular,
we write $G_n$ for $G$.

For a positive integer $n$ let $Q_{n-1}$ denote the following subgroup
of $\GL(n,q)$.
\begin{equation}
\label{Qnminus1}
Q_{n-1} = \left\{
\left[
\begin{array}{cc} 1 & v^t \\ 0 & x \end{array}
\right]
\mid v \in \mathbb{F}_q^{n-1}, x \in \GL(n-1,q) \right\}
\end{equation}
(By convention, $Q_0$ is the trivial subgroup of $\GL_1(q)$.)
Thus $Q_{n-1}$ is the affine group of degree $n-1$. We identify $Q_{n-1}$
with the semidirect product $V_{n-1}G_{n-1}$, where $V_{n-1}$ is the
{\em unipotent radical} of $Q_{n-1}$, consisting of those matrices
in (\ref{Qnminus1}) with $x = 1$.

Suppose now that $n \geq 2$.
Since $G_{n-1}$ acts transitively on the non-identity elements of $V_{n-1}$,
there are two types of irreducible characters of $Q_{n-1}$. The first type
consists of the characters of $G_{n-1}$, inflated to characters of
$V_{n-1}G_{n-1}$. For the second type, we choose a particular element
$\lambda \in \Irr(V_{n-1})$ such that the stabilizer of $\lambda$ in
$G_{n-1}$ equals $Q_{n-2}$. Then the irreducible characters of
$V_{n-1}G_{n-1}$, which do not have $V_{n-1}$ in their kernel, are
parametrized by the irreducible characters of $Q_{n-2}$. We write $\psi_\mu$
for an irreducible character of the second type with parameter $\mu \in
\Irr(Q_{n-2})$. Thus $\psi_\mu = (\hat{\lambda} \cdot
\tilde{\mu})^{V_{n-1}G_{n-1}}$, where $\hat{\lambda}$
is a trivial extension of $\lambda$ to its stabilizer $V_{n-1}Q_{n-2}$, and
$\tilde{\mu} := \infl_{V_{n-1}Q_{n-2}}(\mu)$ is the inflation of $\mu$ to this
stabilizer.

We choose the irreducible character $\lambda$ of $V_{n-1}$ as follows.
Let $U_n$ denote the group of upper triangular unipotent matrices in $G_n$.
Choose a non-trivial homomorphism $\nu : \mathbb{F}_q \rightarrow
\mathbb{C}^*$. Then let $\lambda \in \Irr(U_n)$ be defined by
$\lambda(u) = \prod_{i = 1}^{n-1} \nu(u_{i,i+1})$ for $u = (u_{ij})
\in U_n$. Then $\lambda^{G_n} = \gamma_n$, the character of the
Gelfand Graev representation of $G_n$, (see \cite[Section 8.1]{C}). We also
denote by the same letter the restriction of $\lambda$ to any subgroup of
$U_n$, in particular to the subgroup $V_{n-1}$.

With this notation we are now going to define, recursively on $n-1$
and $i$, $0 \leq i \leq n-1$, the {\em level-$i$-Steinberg character}
$\sigma_i^{(n-1)}$ of $Q_{n-1}$. To begin with, $\sigma_0^{(0)}$ is the
trivial character
of the trivial group $Q_0$. For $n \geq 2$ and $i = 0$, we let
$\sigma_0^{(n-1)}$ denote the inflation of $\St_{n-1}$ to $V_{n-1}G_{n-1}$,
and call it the level-$0$-Steinberg character of $Q_{n-1} = V_{n-1}G_{n-1}$.
For $i \geq 1$, the level-$i$-Steinberg character of $Q_{n-1}$ is defined
by $\sigma_i^{(n-1)} := \psi_\mu$ for $\mu = \sigma_{i-1}^{(n-2)}$.

With this notation we can state our first result. This is a special case of
the results of \cite[Chapter~$5$]{BDK}.
\begin{propo}\label{Propo71}
For all $n \geq 1$, we have
$(\St_{n})_{Q_{n-1}} = \sum_{i=0}^{n-1} \sigma_i^{(n-1)}$.
\end{propo}
\begin{prf}
It is clear, that among the constituents of $(\St_{n})_{Q_{n-1}}$
of the first type, only the inflation of the Steinberg character
$\St_{n-1}$ occurs, and this with multiplicity $1$. The result is
trivial for $n = 1$. Suppose that $n \geq 2$ and let $\mu \in
\Irr(Q_{n-2})$. Using the facts that $V_{n-1}Q_{n-2}G_{n-1} = Q_{n-1}$
and $V_{n-1}Q_{n-2} \cap G_{n-1} = Q_{n-2}$, as well as $(\St_{n})_{Q_{n-1}}
= (\St_{n-1})^{Q_{n-1}}$ (see \cite[Proposition 6.3.3]{C}), we compute
\begin{eqnarray*}
\left( (\St_{n})_{Q_{n-1}}, \psi_\mu \right)
& = & \left( (\St_{n-1})^{Q_{n-1}}, \psi_{\mu} \right) \\
& = & \left( (\St_{n-1})^{Q_{n-1}},
(\hat{\lambda} \cdot \tilde{\mu})^{Q_{n-1}} \right) \\
& = & \left( \St_{n-1},
((\hat{\lambda} \cdot \tilde{\mu})^{Q_{n-1}})_{G_{n-1}} \right) \\
& = & \left( \St_{n-1}, ( (\hat{\lambda} \cdot \tilde{\mu})_{V_{n-1}Q_{n-2} \cap G_{n-1}} )^{G_{n-1}} \right) \\
& = & \left( \St_{n-1}, \mu^{G_{n-1}} \right) \\
& = & \left( (\St_{n-1})_{Q_{n-2}}, \mu \right ).
\end{eqnarray*}
By induction, $(\St_{n-1})_{Q_{n-2}} = \sum_{i=0}^{n-2}
\sigma_i^{(n-2)}$, and the result follows.
\end{prf}

Let $\hat{\omega} := \hat{\omega}_n$ denote the permutation character of $G_n$
on its natural vector space. Thus $\hat{\omega}_n = 1_G + (1_{Q_{n-1}})^G$. 
Hence $\hat{\omega}_n \cdot \St_n = \St_n + ((\St_n)_{Q_{n-1}})^{G_n}$.

Recall that $\gamma_n$ denotes the Gelfand-Graev
character of $G_n$. (For $n = 1$, $\gamma_1$ equals the regular
character of $G_1 = \GL(1,q)$.) For $0 \leq m \leq n$, we let $P_{m}$
denote the standard parabolic subgroup of $G_n$ corresponding to the
composition $(m,n-m)$ of $n$. The unipotent radical of $P_{m}$
is denoted by $U_{m,n-m}$. The Levi subgroup of $P_{m}$ is isomorphic
to $G_m \times G_{n-m}$ (with the convention that $G_0$ denotes the
trivial group).


\begin{theo}\label{ThmGLn}
Let $n \geq 1$ and $G = \GL(n,q)$. Then
$$\hat{\omega} \cdot \St =
\sum_{m=0}^n \left( \infl_{P_{m}}\left(\St_m \boxtimes \gamma_{n-m}\right)\right)^{G}.$$
\end{theo}
\begin{prf}
The summand for $m = n$ on the right hand side equals $\St = \St_n$.
So it suffices to prove that
$$(1_{Q_{n-1}})^G \cdot \St_n = \sum_{m=0}^{n-1}
\left( \infl_{P_m}\left(\St_m \boxtimes \gamma_{n-m}\right)\right)^{G}.$$
Now $(1_{Q_{n-1}})^G \cdot \St_n =
((\St_n)_{Q_{n-1}})^G =
(\sum_{i=0}^{n-1} \sigma_i^{(n-1)})^G$ by Proposition~\ref{Propo71}. To complete
the proof we show that
$$(\sigma_i^{(n-1)})^G =
\left(\infl_{P_{n-i-1}}\left(\St_{n-i-1} \boxtimes \gamma_{i+1}\right)\right)^{G}$$
for all $0 \leq i \leq n - 1$.

Let us start with the case $i = 0$. Here, $Q_{n-1}$ is a normal subgroup
of $P_{1}$, in fact $Q_{n-1} = V_{n-1}G_{n-1}$ and
$P_{1} = V_{n-1}(G_1 \times G_{n-1})$ (in fact $V_{n-1} = U_{1,n-1}$).
Hence $({\sigma_0^{(n-1)}})^{P_{1}} = \infl_{V_{n-1}}( \rho_{G_1} \boxtimes
\St_{n-1})$. It follows that $({\sigma_0^{(n-1)}})^G = \left(\infl_{P_1}\left(
\gamma_1 \boxtimes \St_{n-1}\right)\right)^G$, as claimed.

For $i \geq 1$ (and hence $n \geq 2$) consider the subgroup
$H := (U_{i+1} \times G_{n-i-1})U_{i+1,n-i-1}$
of $P_{i+1} = (G_{i+1} \times G_{n-i-1})U_{i+1,n-i-1}$ (recall
that $U_m$ denotes the group of upper triangular unipotent matrices
in $G_m$). Clearly, $H \leq Q_{n-1}$. We claim that
$\left(\infl_H\left(\lambda \boxtimes \St_{n-i-1}\right)\right)^{Q_{n-1}} =
\sigma_i^{(n-1)}$.

Suppose that this claim has been proved. Then
\begin{eqnarray*}
({\sigma_i^{(n-1)}})^G & = &
\left(\infl_H \left(\lambda \boxtimes \St_{n-i-1}\right)\right)^G \\
 & = & \left(\left(\infl_H\left(
 \lambda \boxtimes \St_{n-i-1}\right)\right)^{P_{i+1}}\right)^G \\
 & = & \left(\infl_{P_{i+1}}\left(\gamma_{i+1} \boxtimes \St_{n-i-1}\right)\right)^G,
\end{eqnarray*}
giving the result.

It suffices to prove the above claim. First observe that $\sigma_i^{(n-1)}(1)
= (q^{n-1} - 1)(q^{n-2} - 1) \cdots (q^{n-i} - 1)\St_{n-i-1}(1)$, and that
this number also equals the degree of the induced character
$\left(\infl_H\left(\lambda \boxtimes \St_{n-i-1}\right)\right)^{Q_{n-1}}$.
By definition, $\sigma_i^{(n-1)} = (\hat{\lambda} \cdot \tilde{\mu})^{Q_{n-1}}$
with $\tilde{\mu} = \infl_{V_{n-1}Q_{n-2}}(\sigma_{i-1}^{(n-2)})$.
Since $i \geq 1$, we have $H \leq V_{n-1}Q_{n-2}$, and thus it suffices to show
that $\hat{\lambda} \cdot \infl_{V_{n-1}Q_{n-2}}({\sigma_{i-1}^{(n-2)}})$ is a
constituent of $(\infl_H(\lambda \boxtimes \St_{n-i-1}))^{V_{n-1}Q_{n-2}}$.
By Frobenius reciprocity, we are left to show that
$\left(\hat{\lambda} \cdot \infl_{V_{n-1}Q_{n-2}}({\sigma_{i-1}^{(n-2)}})\right)_H$
contains $\infl_H(\lambda \boxtimes \St_{n-i-1})$ as a constituent. This is done
by induction on $n$, the case $n = 2$ being trivial.

Since $V_{n-1} \leq H \leq V_{n-1}Q_{n-2}$, we have $H = V_{n-1}K$
with $K = H \cap Q_{n-2}$. Now $H/V_{n-1} \cong K =
(U_i \times G_{n-i-1})U_{i,n-i-1}$,
and, by induction, $\lambda \boxtimes \St_{n-i-1}$ is a constituent of
the restriction of $\sigma_{i-1}^{(n-2)}$ to $K$ (where $\lambda$ is
considered as a character of $U_i$). By the definition of $\lambda$
and of $\hat{\lambda}$ above, it follows that the restriction of
$\hat{\lambda} \cdot \infl_{V_{n-1}Q_{n-2}}({\sigma_{i-1}^{(n-2)}})$
to $H$ contains $\lambda \boxtimes \St_{n-i-1}$ as a constituent.
This completes the proof.
\end{prf}

Multiplying the expression for $\hat{\omega}$ in Theorem~\ref{ThmGLn} by 
$1_G^-$, yields the statement in Theorem~\ref{Thm1Intro} for the general
linear groups.

By this theorem, $\hat{\omega} \cdot \St$ is not multiplicity free,
since every $\gamma_{n-m}$ contains $\St_{n-m}$ as a constituent,
and $\infl_{P_m}\left(\St_{m} \boxtimes \St_{n-m}\right)^G$ contains
$\St_G$ as a constituent. (By \cite[Theorem (70.24)]{CRII}, the latter
assertion can be transformed to a statement in the symmetric group~$S_n$,
where it is obvious.)

\section{Applications}

In this section we prove Theorems~\ref{Thm2Intro} and~\ref{Thm3Intro}.

\subsection{Restricting the Steinberg character}\label{RestrictingSteinberg}
If $G=\Sp(2n,q)$, $q$ odd, and~$P$ denotes the stabilizer of a line in the 
natural module
of~$G$, the characters of $P$ have been described recursively in \cite{AH}.
Rather than recalling the details of~\cite{AH}, we discuss the corresponding
problem for the unitary groups, which reveals a new type of problem. Thus let
$G = U(d,q)$ acting on the vector space $V$ equipped with the Hermitian form
as in Subsection~\ref{UnitaryGroups}. Let $P$ be the stabilizer of an isotropic
line of~$V$. Let $U$ be the unipotent radical of $P$, and let $Z(U)$ denote the
centre of $U$. Additionally, let $L$ denote a Levi subgroup of $P$. Then $P=LU$
and $L = L' \times A$ with $L' \cong U(d-2,q)$ and $A \cong \GL(1,q^2)$.

There are three types of characters of~$P$:
\begin{itemize}
\item[] Type (A): The characters  trivial on $U$.
\item[] Type (B): The characters non-trivial on $U$ but trivial on $Z(U)$.
\item[] Type (C): The characters non-trivial on $Z(U)$.
\end{itemize}

It is slightly less technical to work with the group $P' := L'U$, the
stabilizer of an isotropic vector. Thus $P'$ is a normal subgroup of~$P$ with
cyclic quotient generated by~$A$. This fact can be used to extend the results
below from~$P'$ to~$P$. Of course, the above classification of the irreducible
characters also holds for~$P'$.

Set $\ov{U}:=U/Z(U)$. Observe that $\ov{U}$ is an abelian group which can
be viewed as the natural $\mathbb{F}_{q^2}L'$-module (that is,
$\mathbb{F}^{d-2}_{q^2}$). 
The group $\Irr(\ov{U})$ of irreducible characters of $\ov{U}$ is isomorphic to $\ov{U}$
as abelian groups and as $\mathbb{F}_{q^2}L'$-modules.
In particular, if $\lam \in \Irr(U)$ then the stabilizer of $\lam$ in $L'$
(or the inertia group) coincides with the stabilizer in $L'$ of some element
$\ov{U}$. This simplifies the study of the inertia groups.

Let $\chi$ be an irreducible character of~$P'$ Type $B$. By Clifford's theorem,
there is a non-trivial irreducible character $\lam$ of $\ov{U}$ such that $\chi$ is
induced from an irreducible character $\mu$, say, of the stabilizer $P'_\lam$  of
$\lam$ in $P'$.

So the first matter is to describe $P'_\lam$. It has been observed above
that $P'_\lam =U{\rm Stab}_{L'}(\lam )$ and the second group here coincides
with the stabilizer of some non-zero vector $v \in \mathbb{F}^{d-2}_{q^2}$. As
$\mathbb{F}^{d-2}_{q^2}$ is the natural $\mathbb{F}_{q^2}L'$-module, this
space possesses a unitary form, so the vector in question can be either
isotropic or anisotropic. The group $L' \cong U(d-2,q)$ acts transitively
on the set of (non-zero) isotropic vectors and has $q-1$ orbits on the set
of anisotropic vectors, so~$L'$ has exactly~$q$ orbits on the non-zero vectors
of $\mathbb{F}^{d-2}_{q^2}$. Since the stabilizers of proportional vectors
are the same, we may assume that the representatives of the orbits of
anisotropic vectors all have the same stabilizer.
Thus, $\chi$ corresponds either to an isotropic or to an anisotropic vector.
Depending on this, we say that $\chi$ is of Type $(B1)$ or $(B2)$.

According to this, $P'_\lam/U$ is isomorphic either to $P'_{d-2}$, where
$P'_{d-2}$ is the stabilizer in $L' \cong U(d-2,q)$ of an isotropic vector,
thus defined analogously to~$P'$, or $P'_\lam/U$ is isomorphic to $U(d-3,q)$,
the stabilizer of an anisotropic vector. Thus the irreducible
characters of~$P'$ of Type (B1) are naturally labelled (bijectively)
by $\Irr(P'_{d-2})$, and those of Type (B2) by $\Irr(U(d-3,q))$
(for more details see \cite[2.3.2]{AH}). The characters of Type (B)
are invariant in~$P$, so each of them has exactly $q - 1$ extensions
to~$P$.

Let $\chi = \chi_\mu$ be an irreducible character of~$P$ of Type (B),
labelled by the irreducible character~$\mu$ of $H \leq U(d-2,q)$,
with $H = P'_{d-2}$ or $U(d-3,q)$, respectively. As in \cite[Section~$3$]{AH},
we have
\begin{equation}
\label{TypeB}
(\St_{P'}, \chi_\mu) = (\St^{(L')}_H,\mu).
\end{equation}
We have a similar result as in \cite{AH} for characters of Type (C).
These can be labelled by $\Irr(L')$, such that $\vartheta \in \Irr(L')$
determines exactly $q - 1$ irreducible characters ${\psi}^i_\vartheta$ of Type (C),
permuted transitively by the action of~$P$. Thus every ${\psi}^i_\vartheta$ induces
to an irreducible character~$\psi_\vartheta$ of~$P$ of Type~(C), whose
restriction to~$P'$ equals $\sum_{i=1}^{q-1}{\psi}^i_\vartheta$.

If $\chi$ is an irreducible character of this type labelled by the pair
$(\vartheta,i)$ with $\vartheta \in \Irr(U(d-2,q))$ and $1 \leq i \leq q - 1$,
we have
\begin{equation}
\label{TypeC}
(\St_{P'}, \chi) = (\omega' \cdot\,\St_{L'},\vartheta),
\end{equation}
where $\omega'$ denotes the Weil character of $U(d-2,q)$.
In particular, this multiplicity is independent of~$i$ and
can be computed by Theorem~\ref{Thm1Intro}.

\medskip

\noindent
\textbf{Proof of Theorem~\ref{Thm2Intro}}. If $G = \GL(n,q)$, the result 
follows from \cite[Chapter~$5$]{BDK}. An explicit version is given in
Proposition~\ref{Propo71}.

Next let $G = U(d,q)$.
Suppose first that $H' \leq G$ is the stabilizer of an anisotropic vector.
By the result of Brunat (see the appendix), $\St_{H'}$ equals the product of
the Steinberg character and the Weil character of~$H'$. Using
Corollary~\ref{MultiplicityFree}, the result follows in this case.

Now suppose that $P' \leq G$ is the stabilizer of a non-zero isotropic vector.
Clearly, the only character
of Type~(A) contained in $\St_{P'}$ equals $\St_{L'}$, and it occurs with
multiplicity~$1$.
Now lets look at characters of Type (B). For characters of Type (B1)
we use~(\ref{TypeB}) and induction on~$d$ (the case of $d = 2$ being
clear). For characters of Type (B2) we have to determine the
restriction of the Steinberg character of $L' = U(d-2,q)$ to its
subgroup $L'' = U(d-3,q)$. By what we have proved already, this
restriction is multiplicity free.
The assertion for characters of Type~(C) follows from~(\ref{TypeC})
together with Corollary~\ref{MultiplicityFree}.

Finally, let $G = \Sp(2n,q)$, and let~$P$ denote the stabilizer of a
line $\langle v \rangle$. If $n \leq 3$ the result is already contained
in~\cite{AH}. In the general case it follows from \cite[Corollary~$3.3$]{AH},
together with Corollary~\ref{MultiplicityFree}. Now suppose that $P' \leq P$
is the stabilizer of the vector~$v$. For characters of Type~$3$ (Notation
from~\cite{AH}), the claim easily follows from Clifford theory applied to
the normal subgroup~$P'$ of~$P$. For characters of Type~$2$ we could also
use Clifford theory, but it is simpler to use exactly the same direct
approach as in the unitary groups for characters of Type~(B1).

This completes the proof of Theorem~\ref{Thm2Intro}.

\subsection{The decomposition of a projective character} 

In order to prove Theorem~\ref{Thm3Intro},
we continue our investigation of the generalized spinor representation of the
symplectic groups in characteristic~$2$ begun in Subsection~\ref{ModularCase}.
In particular, we use the notation summarized there. Moreover, we let $st$ 
denote the Brauer reduction modulo~$2$ of the Steinberg representation of
$G = \Sp(2n,q)$, where $q$ is a power of~$2$. Then $st$ is a projective
$\mathbb{F}_qG$-representation. Hence every representation of the form
$st \otimes \phi$ is also projective for every representation~$\phi$ of
$\mathbb{F}_qG$.

\subsubsection{The product $\sigma_n\otimes \sigma_n$ and the natural 
permutation module $\Pi_n$}

\begin{lemma} \label{m76}
The multiplicity of every irreducible $\mathbb{F}_2$-representation $\tau$ of
$\Sp(2n,2)$ in ${(\phi_{\lambda_n}\otimes \phi_{\lambda_n})}_{\Sp(2n,2)}$ is 
equal to the multiplicity of  $\tau$ in the permutation module $\Pi_n$ of 
$\Sp(2n,2)$ associated with the natural action of $\Sp(2n,2)$ on the vectors
of its standard module $V$ (the zero vector is not excluded).
\end{lemma}
\begin{prf}
It suffices to show that the Brauer characters of the two modules 
coincide. 
The action of the image $\eta(\ESp(2n,2))$ of the extrasymplectic group 
on the set of matrices $\mbox{\rm Mat}(2^n,\mathbb{C})$ by conjugation turns
$\mbox{\rm Mat}(2^n,\mathbb{C})$ into a $\mathbb{C}\ESp(2n,2)$-module. Of 
course, this is 
exactly the module afforded by $\overline{\eta} \otimes \eta$. By Corollaries
\ref{pc8} and \ref{c75}, the character of this module at an odd order element
$g \in \ESp(2n,2)$ is equal to $2^{N(V;h)}$ where $h$ is the projection of $g$
into $\Sp(2n,2)$ and $N(V;h)$ is the dimension of the $1$-eigenspace of~$h$ 
on~$V$. Obviously, this coincides with the character of $h$ on $\Pi_n$. By 
Proposition \ref{p88}, the Brauer reduction modulo~$2$ of~$\eta$ equals 
${(\phi_{\lambda_n})}_{\Sp(2n,2)}$. Hence, by Corollary \ref{pc8}, the 
reduction modulo~$2$ of $\overline{\eta} \otimes \eta$ has the same Brauer 
character as ${(\overline{\phi}_{\lambda_n} \otimes 
\phi_{\lambda_n})}_{\Sp(2n,2)}$ and this coincides with ${(\phi_{\lambda_n} 
\otimes \phi_{\lambda_n})}_{\Sp(2n,2)}$ as $\eta$ is real. So the Brauer 
character of ${(\phi_{\lambda_n} \otimes \phi_{\lambda_n})}_{\Sp(2n,2)}$ 
coincides with the Brauer character of the permutation module in question.
\end{prf}

Observe that the natural permutation $\mathbb{F}_q\Sp(2m,2^k)$-module can be
identified with the restriction of $\Pi_{mk}$ to $\Sp(2m,2^k)$, where 
$\Pi_{mk}$ is the natural permutation $\mathbb{F}_2\Sp(2mk,2)$-module.

\begin{lemma} \label{sz9}
For $0\leq i\leq 2n$ let $V_i$ denote the $i$-th exterior power of $V$,
the natural $\mathbb{F}_2\SL(2n,2)$-module ($V_0$ is regarded as 
the trivial module).
Let $\tau$ be an $\mathbb{F}_2\Sp(2n,2)$-composition factor of $\Pi_n$.

$(1)$ Then $\tau$ is isomorphic to a composition factor of ${(V_i)}_{\Sp(2n,2)}$ 
for some $i\leq n$. 

$(2)$ If $\mu$ is a composition factor of $\Pi_n$ viewed as 
$\mathbb{F}_q\Sp(2m,q)$-module, where $q=2^k$ and $n=mk$ then $\mu
= {(\phi_{\lambda})}_{\Sp(2m,q)}$ for $\lambda = \sum_{i=0}^{k-1}
2^i\lambda_{j_i}$ with $j_i\in \{0\ld n\}$. (Recall that the $\lambda_i$ are
the fundamental weights for $i=1\ld n$ and $\lambda_0=0$.)

$(3)$ There is at most one composition factor in $(2)$ occuring with 
multiplicity~$1$; this is ${(\phi_{\lambda})}_{\Sp(2m,q)}$ where $\lambda = 
\sum_{i=0}^{k-1} 2^i\lambda_{m}=(q-1)\lambda_m$. 
\end{lemma}
\begin{prf}
(1) and (2) are proved in \cite[Proposition 3.5]{sz90}. To justify (3), 
consider $\Pi_n$ and $V_i$ as $\SL(2n,2)$-modules, and consider $V_i$ as 
$\mathbb{F}_q\SL(2n,q)$-module. The composition factors of ${(\Pi_n)}_{\SL(2n,2)}$
are irreducible $\mathbb{F}_2\SL(2n,2)$-modules isomorphic to ${(V_i)}_{\SL(2n,2)}$ for 
$i=0\ld 2n-1$, where each factor occurs with multiplicity~$1$ except for the 
trivial one which occurs twice. (This is well known but one may consult 
\cite[Theorem 1.4]{sz90}, where the composition factors of the permutation 
module of $\SL(m,q)$ on the vectors of the natural module have been determined.)
Therefore, the multiplicity of every composition factor in ${(\Pi_n)}_{\Sp(2n,2)}$
and in $\oplus_{i=0}^{2n}{(V_i)}_{\Sp(2n,2)}$ coincide. It is well known that $V_i$
and $V_{2n-i}$ are dual $\SL(2n,2)$-modules. Therefore, ${(V_i)}_{\Sp(2n,2)} \cong 
{(V_{2n-i})}_{\Sp(2n,2)}$. It follows that the irreducible constituents of 
multiplicity~$1$ can only occur in ${(V_n)}_{\Sp(2n,2)}$. Observe that 
${(V_i)}_{\Sp(2n,2)}$ for $i\leq n$ contains a composition factor 
$W_i$ of highest weight $\lambda_i$. By (1) only ${(W_n)}_{\Sp(2n,2)}$ can occur 
in ${(\Pi_n)}_{\Sp(2n,2)}$ with multiplicity~$1$. This completes the case $q=2$.
In general, it follows from this that only irreducible constituents of 
${(W_n)}_{\Sp(2m,q)}$ can occur with multiplicity~$1$. By Lemma \ref{nd5},
${(W_n)}_{\Sp(2m,q)}$ is irreducible and coincides with ${(\phi_{\lambda})}_{\Sp(2m,q)}$
where $\lambda$ is as in Statement~(3).
\end{prf}
\begin{remar}
{\rm
(1)  In fact, the composition factor $\phi_\lambda$ in~(3) occurs with
multiplicity~$1$. This can be proved straightforwardly but we will deduce it
later from Corollary~\ref{MultiplicityFree}.
Observe that Corollary~\ref{c20} implies that 
the composition factors of $\phi_{(q-1)\lambda_n}\otimes \phi_{(q-1)\lambda_n}$
and ${(\Pi_n)}_{\Sp(2m,q)}$ have the same multiplicities.

(2) The composition factors of ${(V_i)}_{\Sp(2n,2)}$ are also studied by Baranov 
and Suprunenko in \cite{BS}.
}
\end{remar}

\subsubsection{Indecomposable summands of $\sigma_n \otimes st$}

In this section we determine the indecomposable constituents of $\sigma_n
\otimes st$.
Let $\nu$ be a dominant weight. We denote by $\phi_\nu$ the irreducible 
representation of $\mathbf{G}$ with highest weight $\nu$. Recall that every 
irreducible representation of $G = \Sp(2n,q)$ is of shape ${(\phi_\nu)}_G$ where
$\nu$ is a $q$-restricted dominant weight of $\mathbf{G} = \Sp(2n,\mathbf{K})$.
Put $\tilde \om := \lambda_1 + \cdots + \lambda_n$.
It is well known that $(q-1)\tilde\om$ is the only $q$-restricted dominant 
weight $\rho$ such that ${(\phi_\rho)}_G=st$. Recall that $\sigma_n = 
{(\phi_{(q-1)\lambda_n})}_G$ and that $\phi_{(q-1)\lambda_n}$ is self-dual. 

\begin{lemma} \label{ch9}
{\rm \cite[9.4]{Hu}}
Let $\psi$ be an irreducible $\mathbb{F}_q G$-module. Then the multiplicity of 
the principal indecomposable module $\Phi_\nu$ in $\psi \otimes st$ is equal
to the multiplicity of $st$ in ${(\phi_\nu)}_G \otimes \psi^*$ where $ \psi^*$
is the dual of $\psi$.
\end{lemma}

There is further information on those $\nu$ for which $\Phi_\nu$ may actually 
occur as a direct summand of $\psi \otimes st$, see \cite[9.4]{Hu}. We could 
prove Theorem~\ref{Thm3Intro} on the base of that information but our special
case can probably be dealt with more efficiently staightforwardly. (Our 
argument here is based on Lemma \ref{sz9} and general facts on representations
of algebraic groups.)

Set $\nu = a_1 \lambda_1 + \cdots + a_n \lambda_n$ where $0 \leq a_1 \ld
a_n \leq q-1$, and  $\nu' = a_1 \lambda_1 + \cdots + a_{n-1} \lambda_{n-1}$.

\medskip
{\bf Proof of Theorem \ref{Thm3Intro}.} We show that  $\Phi_\nu$ is a direct
summand of $\sigma_n \otimes st$ if and only if  
$\nu' = (q-1)(\lambda_1 + \cdots + \lambda_{n-1})$, that is, $a_1 = \cdots =
a_{n-1} = q-1$. It can be deduced from Steinberg \cite[Corollary to Theorem
41 and Theorem 43]{St} that $\phi_{\nu'} \otimes \phi_{(q-1)\lambda_n} = 
\phi_{\nu' + (q-1)\lambda_n}$. If $a_n=0$, we have $\nu = \nu'$ so the 
representation $\phi_{\nu + (q-1)\lambda_n}$ is irreducible. As $\nu +
(q-1)\lambda_n$ is a dominant $q$-restricted weight, 
${(\phi_{\nu + (q-1)\lambda_n})}_G$ is irreducible, so it is not equal to $st$
unless $\nu = (q-1)(\lambda_1 + \cdots + \lambda_{n-1})$. So the claim follows
from Lemma \ref{ch9}.

Next assume $a_n>0$. Then we have that 
$$\phi_{\nu} \otimes \phi_{(q-1)\lambda_n} = \phi_{\nu'} 
\otimes \phi_{a_n\lambda_n} \otimes \phi_{(q-1)\lambda_n}.$$ 
Let $a_n = \sum_{i=0}^{k-1} 2^ib_i$ be the $2$-adic expansion of $a_n$ (so 
$0 \leq b_i \leq 1$). Then 
$$\phi_{a_n\lambda_n} \otimes \phi_{(q-1)\lambda_n}=(\phi_{b_0\lambda_n} \otimes
\phi_{\lambda_n}) \otimes F_0(\phi_{b_1\lambda_n} \otimes \phi_{\lambda_n})
\otimes \cdots \otimes F_0^{k-1} (\phi_{b_{k-1}\lambda_n} \otimes 
\phi_{\lambda_n}).$$ 
If $b_i = 0$ then $\phi_{b_{i}\lambda_n} \otimes \phi_{\lambda_n} = 
\phi_{\lambda_n}$, otherwise $b_i = 1$ and the composition factors of 
${(\phi_{b_{i}\lambda_n} \otimes \phi_{\lambda_n})}_G$ are ${(\phi_{\lambda_j})}_G$
for $0 \leq j \leq n$ by Lemma \ref{sz9}. Therefore, the composition factors of
${(\phi_{a_n\lambda_n} \otimes \phi_{(q-1)\lambda_n})}_G$ are the restrictions 
to~$G$ of representations of shape
$$\phi_{\lambda_{i_0}} \otimes F_0(\phi_{\lambda_{i_1}}) \otimes \cdots \otimes
F_0^{k-1}(\phi_{\lambda_{i_{k-1}}}) = \phi_{\lambda_{i_0} + 2\lambda_{i_1} + \cdots
 + 2^{k-1} \lambda_{i_{k-1}}}$$ 
where $0 \leq i_0, i_1, \ldots ,i_{k-1} \leq n$.
Moreover,  Lemma \ref{sz9} tells us that the multiplicity of 
${(\phi_{\lambda_j})}_G$ in ${(\phi_{b_{i}\lambda_n} \otimes \phi_{\lambda_n})}_G$ 
(when $b_i = 1$) is at least~$2$ unless $j = n$. Therefore every composition 
factor $\tau$, say, of 
$${(\phi_{\nu'} \otimes \phi_{\lambda_{i_0} + 2\lambda_{i_1} + \cdots + 2^{k-1} 
\lambda_{i_{k-1}}})}_G$$ 
occurs at least twice unless $\lambda_{i_0} = \lambda_{i_1} = \cdots = 
\lambda_{i_{k-1}} = \lambda_n$ in which case $\lambda_{i_0} + 2\lambda_{i_1} +
\cdots + 2^{k-1} \lambda_{i_{k-1}} = (q-1)\lambda_n$. It follows that $\tau 
\neq st$ if $\tau$ occurs more than once, 
as otherwise, by Lemma \ref{ch9}, $\Phi_\nu$ occurs at least twice in 
${(\phi_{\nu} \otimes \phi_{(q-1)\lambda_n})}_G$ which contradicts 
Corollary~\ref{MultiplicityFree}.

So we are left with determining the multiplicity of $st$ in ${(\phi_{\nu'} 
\otimes \phi_{(q-1)\lambda_n})}_G$. As mentioned above, the latter 
representation coincides with ${(\phi_{\nu' + (q-1)\lambda_n})}_G$, which is 
irreducible. It coincides with $st$ if and only if $\nu' = (q-1)(\lambda_1 + 
\cdots + \lambda_{n-1})$.

\begin{remar}
{\rm
The above reasoning justifies also the claim in Remark (1) after Lemma~\ref{sz9}.
}
\end{remar}

\section*{Acknowledgements}

The second author greatfully acknowledges financial support by the
DFG Research Training Group (Gra\-du\-iertenkolleg) ``Hierarchie
und Symmetrie in mathematischen Modellen'', and, at the final stage
of the work, by a Leverhulme Emeritus Fellowship (Grant EM/2006/0030). 

A part of this
work was done during a visit of the first
author at the ``Centre Interfacultaire Bernoulli'' within the
program ``Group Representation Theory'' (January to June 2005).

We thank Frank L\"ubeck for reassuring computations with
{\sf CHEVIE} \cite{chevie} in an early state of this work, as well as
for his careful reading of the manuscript. We also thank Frank
Himstedt for his hint to reference~\cite{AHu}.
Finally we are indebted to Oliver Brunat for pointing out an inaccuracy
in an earlier version of this article.

\end{document}